\documentclass[mathpazo]{cicp}

\usepackage[margin=1.3in]{geometry}
\usepackage{amssymb,amsmath}
\usepackage{amsfonts}
\usepackage{amsthm}
\usepackage{graphicx}
\usepackage{hyperref, xcolor}
\usepackage{mathrsfs}
\usepackage{bm}
\usepackage{subfigure}
\usepackage{float}
\usepackage{indentfirst}
\usepackage{algorithm}
\usepackage{algorithmic}
\usepackage{comment}
\usepackage{appendix}

\begin{document}

\title{\textbf{A threshold dislocation dynamics method}}

\author[]{Xiaoxue Qin\affil{1,4}, Alfonso H.W. Ngan\affil{5}, Yang Xiang\affil{2,3}\corrauth}
\address{\affilnum{1}\  Department of Mathematics, Shanghai University, Shanghai 200444, China\\
\affilnum{2}\ Department of Mathematics, The Hong Kong University of Science and Technology, Clearwater Bay, Kowloon, Hong Kong, China\\
\affilnum{3}\ HKUST Shenzhen-Hong Kong Collaborative Innovation Research Institute, Futian, Shenzhen, China.\\
\affilnum{4}\  Newtouch Center for Mathematics of Shanghai University, Shanghai 200444, China\\
\affilnum{5}\ Department of Mechanical Engineering, The University of Hong Kong, Pokfulam Road, Hong
Kong, China}
\emails{{\tt maxiang@ust.hk} (Y.~Xiang) }

\begin{abstract}
The Merriman-Bence-Osher threshold dynamics method is an efficient algorithm to simulate the motion by mean curvature. It has the advantages of being easy to implement and with high efficiency. In this paper, we propose a threshold dynamics method for dislocation dynamics in a slip plane, in which the spatial operator is essentially an anisotropic fractional Laplacian.  We show that this threshold dislocation dynamics method is able to give { two correct leading orders} in dislocation velocity,  including both the $O(\log\varepsilon)$ local curvature force and the $O(1)$ nonlocal force due to the long-range stress field generated by the dislocations as well as the force due to the applied stress, where $\varepsilon$ is the dislocation core size, { if the time step is set to be $\Delta t=\varepsilon$.
This generalizes the available  result of  threshold dynamics with the corresponding fractional Laplacian, which is on the leading order $O(\log\Delta t)$ local curvature velocity under the isotropic kernel.} We also propose a numerical method based on spatial variable stretching to  correct the mobility and to rescale the velocity for efficient and accurate simulations, which can be applied generally to any threshold dynamics method. We validate the proposed threshold dislocation dynamics method by numerical simulations of various motions and interaction of dislocations.

\end{abstract}

\maketitle

\section{introduction}
Mean curvature flow describes the motion of a co-dimension one object normal to itself with velocity equal to its mean curvature. Merriman, Bence and Osher (MBO) developed an efficient threshold dynamics method to simulate the motion by mean curvature  \cite{merriman1992diffusion,merriman1994motion}. In this method, two simple steps alternates: a convolution with diffusion kernel and a thresholding step. The MBO threshold method has the advantages of being easy to implement and with high efficiency.  The MBO threshold dynamics method has been further developed with some efficient implementations and generalization to multiphase interfaces \cite{ruuth1998efficient,ruuth1998diffusion,esedog2010diffusion,svadlenka2014variational} as well as convergence analysis \cite{evans1993convergence,Barles1995,Laux2016,Yip2017}.
Esedoglu and Otto developed a threshold dynamics method for  dynamics of networks with arbitrary surface tensions \cite{esedog2015threshold}.
Elsey and Esedoglu generalized the threshold dynamics method to anisotropic mean curvature motion by replacing
the isotropic Gaussian kernel in the convolution step of the original algorithm with more
general, anisotropic kernels  \cite{esedoglu2018anisotropy}.
Convergence of nonlocal threshold dynamics corresponding to the fractional Laplacian was proved by Caffarelli and  Souganidis in \cite{caffarelli2010convergence}.
The threshold dynamics method was also extended to { Willmore flow} and some high-order geometric flow problems \cite{esedoglu2008threshold,hu2023unconditionally}, wetting of fluids on rough
surfaces \cite{xu2017efficient},  image segmentation \cite{wang2017efficient}, topology optimization for fluids \cite{WangXP2022,hu2022efficient}, and reconstructing surface from point clouds \cite{wang2021efficient}, etc.

Dislocation dynamics simulation is an important tool for the study of plastic deformation in crystalline materials 
\cite{kubin1992,Ghoniem2000,xiang2003level,quek2006level,cai2006non,Arsenlis2007}, in which the motion and interaction of dislocations (line defects) are simulated. The driving force on dislocations is nonlocal, which is due to the stress field generated by all the dislocations. This is unlike the motion by curvature, which depends only on the local profile of the curve. The driving force on dislocations consists of both the nonlocal $O(1)$ force and the local $O(\log\varepsilon)$ curvature force, where $\varepsilon$ is the dislocation core size, and both are important in the dynamics of dislocations. This force on dislocations is in general anisotropic depending on the orientation of the dislocations.
 The Peierls-Nabarro model \cite{peierls1940size,nabarro1947dislocations,vitek1968intrinsic,anderson2017theory} is a hybrid model that incorporates atomic-size dislocation core  into the continuum framework.  Computational models for dislocation structure and dynamics of curved dislocations based on the Peierls-Nabarro models and generalizations have been developed \cite{Xu2000,Movchan2003,Shen2004,Xiang2006,xiang2008generalized,wei2008generalized,zhu-adaptive,Ngan2018}.
 Velocities of straight dislocations with the applied stress under the setting of fractional Laplacian equation (a simplified Peierls-Nabarro model)
  have been analyzed \cite{Monneau2012,Dipierro2015}.

In this paper, we propose an efficient threshold dynamics method for dislocation dynamics in a slip plane, based on the Peierls-Nabarro model for curved dislocations in \cite{xiang2008generalized,wei2008generalized}. In the convolution step, the dislocation stress field kernel, which is essentially an anisotropic fractional Laplacian kernel,
is used instead of the isotropic diffusion kernel in a standard threshold dynamics method. We show that this proposed threshold dislocation dynamics method gives correct dislocation velocity as compared with the discrete dislocation dynamics method. More precisely, we show that the threshold dislocation dynamics method gives both the correct $O(\log \varepsilon)$ local curvature force and the correct $O(1)$ nonlocal force due to the long-range stress field generated by the dislocations as well as the force due to the applied stress, where $\varepsilon$ is the dislocation core size, { if we set the time step $\Delta t=\varepsilon$ in the threshold dislocation dynamics method.

  Note that based on the diffusion kernel, the threshold dynamics method has been extended to  Willmore flow of planar curves by expansions of two orders of velocities in terms of $\Delta t$ \cite{esedoglu2008threshold,hu2023unconditionally}.
  The velocity of the Willmore flow as well as that of the curvature flow at the leading order in these expansions generated by the diffusion kernel, are both local.  Whereas  the dislocation velocity generated by an anisotropic fractional Laplacian kernel is nonlocal, which is an integral over the entire dislocation.
 In the available  threshold dynamics method with  the corresponding square root Laplacian kernel in \cite{caffarelli2010convergence},
 only the leading order $O(\log \Delta t)$ local curvature velocity under the isotropic kernel was obtained. There is  no analysis for the convergence to the nonlocal velocity under the framework of threshold dynamics methods available in the literature.
Compared with the exponentially-decaying diffusion kernel that generates local velocities, the anisotropic dislocation kernel decays much slower, with $r^{-3}$ where $r$ is the distance between two points, leading  to nonlocal dislocation velocity and making the convergence analysis challenging.
 }

 The proposed threshold dislocation dynamics method can be considered as an efficient implementation of the level set dislocation dynamics method \cite{xiang2003level,quek2006level}, with simple reinitialization and evolution of the level set function, and has the advantages of handling topological and geometrical changes automatically as compared with the front-tracking based discrete dislocation dynamics method, e.g.~\cite{kubin1992,Ghoniem2000,Arsenlis2007}. Compared with the Peierls-Nabarro type models, the threshold dislocation dynamics method is able to use much larger time step ($\Delta t\gg \Delta x$ $vs$ $\Delta t=O(\Delta x)$ or $\Delta t=O(\Delta x^2)$ in Peierls-Nabarro type models and generalizations).

We also develop a variable stretching method to correct the dislocation mobility and to rescale the dislocation velocity in our threshold dislocation dynamics method.
In the threshold dislocation dynamics method, the physical settings for the dislocation velocity to be accurate for the two leading orders impose  restrictions on the numerical implementation of this method.
First, the dislocation stress field kernel, which is in the form an anisotropic fractional Laplacian kernel, leads to an anisotropic dislocation mobility. This anisotropic mobility is not necessarily the dislocation mobility from the physics. Second, the time step $\Delta t$, in the dimensionless form of the equation, reflects the dislocation core radius, which has to be small due to the fact that the dislocation core size should be much less than the size of the domain. As a result, the spatial grid constant $\Delta x$ has to be even smaller due to the requirement $\Delta t\gg \Delta x$ in a threshold dynamics method.
We propose to correct the dislocation mobility by stretching of the spatial variables on the numerical grid, which is different from the available method in the literature \cite{esedoglu2018anisotropy} based on adjustment of the diffusion kernel. The dislocation velocity can also be rescaled to a larger value by this variable stretching method. This variable stretching method can be applied generally to any threshold dynamics method to adjust the velocity and mobility of the moving front.

We perform numerical simulations using the threshold dislocation dynamics method. Simulation results agree with those of theoretic predictions \cite{anderson2017theory} and discrete dislocation dynamics simulations \cite{xiang2003level,xiang2004level}. Especially, these simulation results demonstrate that our threshold dislocation dynamics method can indeed correctly capture both the leading order ($O(\log\varepsilon)$) curvature motion and the next order ($O(1)$) long-range interaction for the dynamics of dislocations.

This paper is organized as follows. We first review the Peierls-Nabarro model in Sec.~\ref{sec:PN}. In Sec.~\ref{sec:Thresh}, we present the threshold dislocation dynamics method, based on the Peierls-Nabarro model, and examine  properties of the dislocation stress kernel in the dislocation dynamics equation. In Sec.~\ref{sec:Ana}, we analyze the dislocation velocity given by the threshold dislocation dynamics method, including  both the velocity due to local curvature on $O(\log \varepsilon)$ and the nonlocal velocity on $O(1)$ due to long-range dislocation interaction.
 In Sec.~\ref{sec:CTDD},
 we present a numerical method based on spatial variable stretching to correct the mobility of dislocations and to speed up the dislocation motion in the threshold dislocation dynamics method.
 The algorithm for the threshold dislocation dynamics method with correction of dislocation mobility and velocity rescaling is presented in Sec.~\ref{sec:algorithm}.
 In Sec.~\ref{sec:Num}, we perform numerical simulations using the threshold dislocation dynamics method
for the motion of a straight dislocation under applied stress, shrinking and expanding of dislocation loops, dislocations bypassing particles, and operation of a Frank-Read source.
Simulation results are compared with  theoretic predictions \cite{anderson2017theory} and discrete dislocation dynamics simulation results \cite{xiang2003level,xiang2004level}.

\section{Review of the Peierls-Nabarro model} \label{sec:PN}
In this section, we review the generalized Peierls-Nabarro model for curved dislocations proposed in Ref.~\cite{xiang2008generalized,wei2008generalized}, as a dynamics model of dislocations.

Suppose that the slip plane of the dislocations is located at $z=0$. We focus on the dislocations with Burgers vector $\bm b=(b,0)$, where $b$ is the magnitude of the Burgers vector. The dislocations are described by the disregistry $\phi(x,y)$ in the direction of the Burgers vector, whose sharp transition regions between regions with values of integer multiples of $b$ represent the dislocation core. The dislocation is in the direction of { $\pmb \xi=\nabla\phi\times\hat{\mathbf z}$, where $\hat{\mathbf z}$} is the unit vector in the $+z$ direction.

The total energy in the framework of the Peierls-Nabarro models \cite{peierls1940size,nabarro1947dislocations,vitek1968intrinsic,hirth1983theory} is:
\begin{equation}\label{eqn:toen}
E=E_{\text{elastic}} +E_{\text{misfit}}.
\end{equation}
The elastic energy $E_{\text{elastic}}$ is:
\begin{equation}
E_{\text{elastic}} =\iiint_{\mathbf R^3}\sum_{i,j=1}^3\frac{1}{2}\sigma_{ij}\epsilon_{ij}dxdydz,
\end{equation}
where $\{\sigma_{ij}\}$ and $\{\epsilon_{ij}\}$ are the stress and strain tensors determined by the disregistry $\phi$. The misfit energy $E_{\text{misfit}}$ due to the nonlinear atomic interaction across the slip plane \cite{vitek1968intrinsic} is:
\begin{equation}
E_{\text{misfit}} =\iint_{\mathbf R^2}\gamma(\phi(x,y))dxdy,
\end{equation}
where $\gamma(\phi)$ is the generalized stacking fault energy, and here we use the Frenkel sinusoidal potential \cite{peierls1940size,nabarro1947dislocations,hirth1983theory}
 \begin{equation}
 \gamma(\phi)=\frac{\mu b^2}{4\pi^2d}\left(1-\cos\frac{2\pi\phi}{b}\right),
 \end{equation}
 where $d$ is the lattice spacing  perpendicular to the slip plane, and $\mu$ is the shear modulus.

The dynamics of dislocations based on the Peierls-Nabarro model is given by the gradient flow of the total energy $E$ in Eq.~\eqref{eqn:toen}, which  is:
\begin{equation}\label{eqn:gradient-f}
\phi_t=-M_p\frac{\delta E}{\delta \phi}=-M_p\left(\sigma_{13}+\frac{\partial\gamma}{\partial\phi}\right),
\end{equation}
where the stress component $\sigma_{13}=\sigma_{13}^{\rm dis}+\sigma_{13}^{\rm app}$, $\sigma_{13}^{\rm dis}$ is the stress  generated by the dislocations:
\begin{flalign}\label{eqn:sigma13-0}
\sigma_{13}^{\rm dis}(x,y)=&\iint_{\mathbf R^2}\left[\frac{\mu }{4\pi(1-\nu)}\frac{(x-\bar{x})}{[(x-\bar{x})^2+(y-\bar{y})^2]^{\frac{3}{2}}} \phi_x(\bar{x},\bar{y})\right.\nonumber\\
& \ \ \ \left.+\frac{\mu }{4\pi}\frac{(y-\bar{y})}{[(x-\bar{x})^2+(y-\bar{y})^2]^{\frac{3}{2}}} \phi_y(\bar{x},\bar{y})\right]d\bar{x}d\bar{y},
\end{flalign}
with $\mu$ being the shear modulus and $\nu$ the Poisson ratio, $\sigma_{13}^{\rm app}$ is the applied stress, and
 $M_p>0$ is the mobility. Note that the dislocation core size is of the order of $b$ or $d$, both of which are of the size of the lattice constant of the crystal.

In the dimensionless form, we use the dimensionless disregistry  $u=\phi/b$, and stretch the length by the length unit of the simulation domain $l_0$, the time by $1/M_p\mu$, and the stress by $\mu b/l_0$. Defining $\varepsilon_0=d/l_0$, the evolution equation \eqref{eqn:gradient-f} for dislocations becomes
\begin{eqnarray}\label{eqn:moelv}
u_t= L(u)-\frac{1}{2\pi \varepsilon_0}\mathrm{sin}(2\pi u)-{\sigma^{\rm app}},
\end{eqnarray}
where $L(u)$ is the dimensionless form of $\sigma_{13}^{dis}$:
\begin{flalign}\label{eqn:sigma13}
L(u)\big|_{(x,y)}=&-\iint_{\mathbf R^2}\left[\frac{ 1}{4\pi(1-\nu)}\frac{(x-\bar{x})}{[(x-\bar{x})^2+(y-\bar{y})^2]^{\frac{3}{2}}} u_x(\bar{x},\bar{y})\right.\nonumber\\
& \ \ \ \left.+\frac{ 1}{4\pi}\frac{(y-\bar{y})}{[(x-\bar{x})^2+(y-\bar{y})^2]^{\frac{3}{2}}} u_y(\bar{x},\bar{y})\right]d\bar{x}d\bar{y},
\end{flalign}
and $\sigma^{\rm app}$ is simplified notation for  $\sigma^{\rm app}_{13}$.
The Fourier transform of $L(u)$ is:
\begin{equation}\label{eqn:sigma13-hat}
\widehat{L(u)}=-\frac{1}{2}\left(\frac{k_1^2}{(1-\nu)\|\mathbf{k}\|} +\frac{k_2^2}{\|\mathbf{k}\|}\right)\hat{u},
\end{equation}
where $\mathbf k=(k_1,k_2)$ is the frequency vector, and $\|\mathbf{k}\|= \sqrt{k_1^2+k_2^2}$. 
Note that the dislocation core size is $O(\varepsilon_0)$, and $\varepsilon_0\ll 1$ in the simulation.

Formulations of the Peierls-Nabarro model for more general cases are reviewed in Appendix~\ref{sec:PN-general}.

\section{Threshold dislocation dynamics method}\label{sec:Thresh}


In this section, we present the threshold dislocation dynamics method. Note that although it
 can be considered as an implementation scheme of the evolution equation \eqref{eqn:moelv} in its form, the purpose of this method is to simulation the dynamics of dislocations, and it will be validated in the next section { by comparing the generated dislocation velocity} with that given by available discrete dislocation dynamics methods~\cite{kubin1992,Ghoniem2000,xiang2003level,Arsenlis2007}. This threshold dynamics formulation will give a specific anisotropic mobility for dislocations. Further correction of the dislocation mobility and {rescaling of the dislocation velocity} will be presented in Sec.~\ref{sec:CTDD}, and the algorithm of this threshold dislocation dynamics method will be summarized in Sec.~\ref{sec:algorithm}.

We first present the method based on a single dislocation $\Gamma$ that evolves in its slip plane $xy$. As in a phase field model, the dislocation and its dynamics are described by the evolution of a function $u(x,y,t)$ over the entire domain. Before and after the evolution of each time step, $u(x,y,t)$ is the characteristic function of the region enclosed by the dislocation $\Gamma$. The dislocation is in the direction of { $\pmb \xi=\nabla u\times\hat{\mathbf z}$, where $\hat{\mathbf z}$ }is the unit vector in the $+z$ direction.

In the first step of the evolution from time $t_n$ to $t_{n+1}=t_n+\Delta t$,
we evolve $u(x,y,t)$ by the following equation
\begin{eqnarray}
&&{u}_t=L({u})-\sigma^{\rm app}, \label{eqn:1-1}\\
&&u|_{t=t_n}=u^n=1_{S_n},\label{eqn:1-2}
\end{eqnarray}
where $1_{S_n}$ is the characteristic function of the region $S_n$ enclosed by the dislocation $\Gamma$ (in the right-hand sense) at time $t_n$.
The solution of Eqs.~\eqref{eqn:1-1} and \eqref{eqn:1-2} at time $t_{n+1}$ is
\begin{equation}\label{eqn:formulation}
{u}(x,y,t_{n+1})= K_{\Delta t}*1_{S_n}-\sigma^{\rm app}\Delta t,
\end{equation}
where $K_{\Delta t}$ is the kernel:
\begin{equation}\label{eqn:K}
K_{\Delta t}(x,y)=\frac{1}{(2\pi)^2}\iint_{\mathbf R^2} e^{-\frac{\Delta t }{2} \left(\frac{k_1^2}{(1-\nu)\|\mathbf{k}\|} +\frac{k_2^2}{\|\mathbf{k}\|}\right)}e^{i(k_1x+k_2y)} \mathrm{d}k_1\mathrm{d}k_2,
\end{equation}
 and $*$ is the convolution operator: 
 $f*g(x,y)=\iint_{\mathbf R^2}f(x-\bar{x}, y-\bar{y})g(\bar{x},\bar{y})\mathrm{d}\bar{x}\mathrm{d}\bar{y}$. (Note that $\widehat{f*g}=(2\pi)^2\hat{f}\hat{g}$ in two dimensions.)
 In a threshold dynamics method, this step is to evolve the function by convolution based on the linear part of the phase field equation~\cite{merriman1992diffusion}, which is Eq.~\eqref{eqn:moelv} in the case of the Peierls-Nabarro model for dislocation dynamics.

{ Note that $\Delta t$ is the core parameter in the kernel  $K_{\Delta t}$ that smooths the characteristic function $1_{S_n}$  of the region $S_n$ enclosed by the dislocation. The  physical meaning of $\Delta t$ is the dislocation core parameter in the dimensionless form presented at the end of Sec.~\ref{sec:PN}. }


In the second step of the evolution from time $t_n$ to $t_{n+1}$, we perform  {thresholding} on the solution $u(x,y,t_{n+1})$ obtained in the first step:
\begin{equation}\label{eqn:thresholding1}
u^{n+1}(x,y)= \left\{
\begin{array}{ll}
0,&{\rm if}\ u(x,y,t_{n+1})\leq 0.5,\vspace{1ex}\\
1,& {\rm if}\ u(x,y,t_{n+1})> 0.5.
\end{array}
\right.
\end{equation}
In a threshold dynamics method, this step is to evolve the function based on the nonlinear part of the phase field equation~\cite{merriman1992diffusion}, which is Eq.~\eqref{eqn:moelv} in the case of  the Peierls-Nabarro model for dislocation dynamics,  in the limit when the core width parameter {$\varepsilon_0$} goes to $0$. Note that $u^{n+1}(x,y)$ is in fact the characteristic function $1_{S_{n+1}}$ of the region $S_{n+1}$ enclosed by the dislocation $\Gamma$  at time $t_{n+1}$, from which the evolution procedure further continues.

For multiple dislocations, especially in the case when one dislocation loop is enclosing another, we can use an integer-valued function, instead of the characteristic function, to represent the dislocations. That is, the convolution step remains the same, and in the { thresholding} step, we have
\begin{equation}\label{eqn:thresholding2}
u^{n+1}(x,y)=
j, \ \ {\rm if}\ j-0.5<u(x,y,t_{n+1})\leq j+0.5,
\end{equation}
for an integer $j$. In this case, the positions of dislocations are represented by the locations of jumps (contour lines of $ {u}=j+0.5 $ before the { thresholding}).

Unlike those available threshold dynamics methods reviewed in the introduction section, all of which focus on the leading order velocity of the moving front that is proportional to its local  curvature on the order of $\log\varepsilon$, here in dislocation dynamics, the velocity is long-ranged and includes contributions from both the leading order $\log\varepsilon$ and the next order $O(1)$
\cite{anderson2017theory,GB1976,Zhao2012}. We will show by asymptotic analysis that the above threshold dislocation dynamics method can indeed generate the correct dislocation velocity in Sec.~\ref{sec:Ana}. Moreover, these two steps of the threshold dislocation dynamics method will lead to a specific anisotropic mobility, and we will present in Sec.~\ref{sec:CTDD} a method by numerical stretching of the spatial variables to adjust the mobility to the desired one (e.g., the isotropic one) and to rescale the velocity to a larger value.
We would like to remark that the method and the analysis can also apply to the more general models discussed in Appendix~\ref{sec:PN-general}.

We summarize this basic threshold dislocation dynamics algorithm below. Algorithm for the method with correction of dislocation mobility and rescaling of dislocation velocity will be presented in Sec.~\ref{sec:algorithm}.

\vspace{0.1in}
\noindent
\underline{\bf Threshold Dislocation Dynamics Method: Basic Algorithm (TDDMB)}
\vspace{0.05in}

1. Give the initial condition $u^0$, and set the time step $\Delta t$  which corresponds to the dislocation core radius.

2. Evolve the solution ${u}$ from $t_n$ to $t_{n+1}$:
\begin{equation*}
\overline{u^{n+1}}= K_{\Delta t}*1_{S_n}-\sigma^{\rm app}\Delta t.
\end{equation*}

3. Update the solution $u$ at $t_{n+1}$ using threshold:
\begin{equation*}
u^{n+1}=j, \ \ {\rm if}\ j-0.5<\overline{u^{n+1}}\leq j+0.5.
\end{equation*}

4. Repeat steps 2-3.

\vspace{0.2in}

Finally in this section, we discuss properties of the kernel function
$K_{\Delta t}(x,y)$ defined in Eq.~\eqref{eqn:K}, which will be used in the analysis of the dislocation velocity in the next section. We first write the kernel function in a general form:
\begin{equation}\label{eqn:Kt}
K(x,y,t)=\frac{1}{(2\pi)^2}\iint_{\mathbf R^2} e^{-\frac{t}{2} \left(\frac{k_1^2}{(1-\nu)\|\mathbf{k}\|} +\frac{k_2^2}{\|\mathbf{k}\|}\right)}e^{i(k_1x+k_2y)} \mathrm{d}k_1\mathrm{d}k_2,
\end{equation}
and especially, $K_{\Delta t}=K(x,y,\Delta t)$. The Fourier transform of $K(x,y,t)$ is $$\widehat{K}(k_1,k_2,t)=\frac{1}{(2\pi)^2}e^{-\frac{t}{2} \left(\frac{k_1^2}{(1-\nu)\|\mathbf{k}\|} +\frac{k_2^2}{\|\mathbf{k}\|}\right)}.$$

\begin{figure}[!htb]
		\centering
		\includegraphics[width=3.2in]{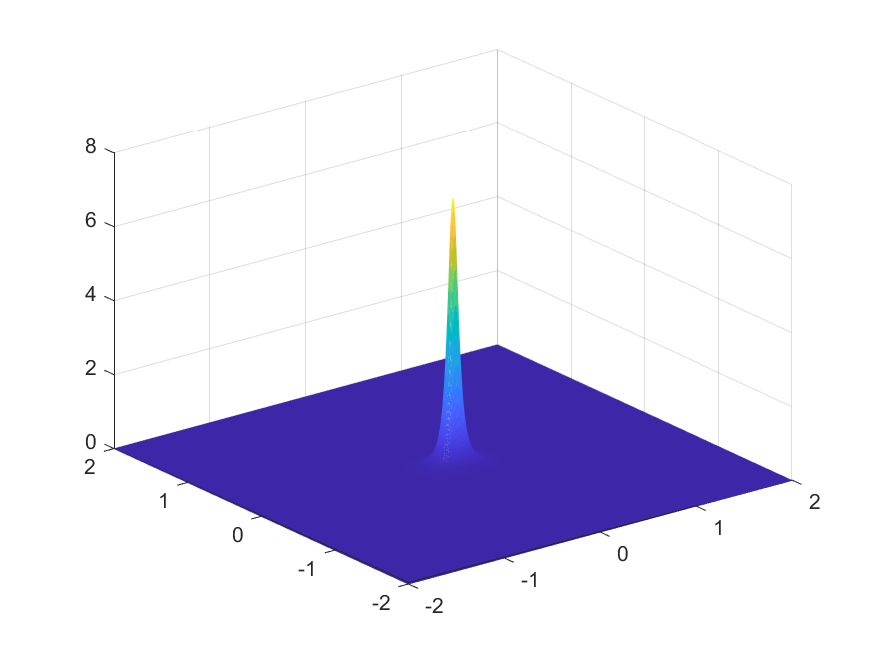}
		\caption{The kernel function $K(x,y,t)$,  with $t=0.1$ and $\nu=0.3$.}
		\label{fig:curve2}
	\end{figure}

The kernel function $K(x,y,t)$ can be considered as a regularized delta-function in two dimensions, with regularization width of $O(t)$. In fact, it is easy to see that $\iint_{\mathbf R^2}K(x,y,t)dxdy=1$, and ${\displaystyle \lim_{t\to 0^+}}K(x,y,t)=\delta (x,y)$, where $\delta (x,y)$ is the two dimensional Dirac delta-function. See Fig.~\ref{fig:curve2} for a plot of $K(x,y,t)$.

When $\nu=0$, $K(x,y,t)$ has the following analytical formula:
\begin{equation}\label{eqn:nu-0}
		K(x,y,t)=K_0(x,y,t)\equiv\frac{1}{(2\pi)^2}\frac{t/2}{[(t/2)^2+x^2+y^2 ]^{\frac{3}{2}}}.
	\end{equation}

When $\nu\neq 0$, analytical expression of $K(x,y,t)$ is not available, and by Eqs.~\eqref{eqn:Kt} and \eqref{eqn:nu-0}, it can be shown that
\begin{equation}\label{eqn:asymp-K}
|K(x,y,t)|\leq\frac{Ct}{(x^2+y^2)^{\frac{3}{2}}},
\end{equation}
for some constant $C$. In fact, letting $G(k_1,k_2)= \frac{1}{(2\pi)^2}e^{-\frac{t}{2} \left(\frac{k_1^2}{(1-\nu)\|\mathbf{k}\|} +\frac{k_2^2}{\|\mathbf{k}\|}\right)}$, from Eq.~\eqref{eqn:Kt}, we have
\begin{flalign*}
\iint_{\mathbf R^2} G(k_1,k_2)e^{i(k_1x+k_2y)}\mathrm{d}k_1\mathrm{d}k_2=&-\iint_{\mathbf R^2} \frac{1}{2}\left(\frac{1}{(ix)^3}\frac{\partial^3 G(k_1,k_2)}{\partial k_1^3}+\frac{1}{(iy)^3}\frac{\partial^3 G(k_1,k_2)}{\partial k_2^3}\right)\\
&\hspace{0.5in}\cdot e^{i(k_1x+k_2y)}\mathrm{d}k_1\mathrm{d}k_2.
\end{flalign*}
  It can be calculated that the  integral on the right-hand side of this equation is integrable and goes to $0$ as $x$ and $y$ go to $+\infty$, and Eq.~\eqref{eqn:asymp-K} follows. The factor $t$ in Eq.~\eqref{eqn:asymp-K} can be obtained by stretching integration variables $k_1$ and $k_2$ in Eq.~\eqref{eqn:Kt}.

 Similarly, it can be shown that
  \begin{equation}\label{eqn:asymp-K1}
\left|\frac{\partial K(x,y,t)}{\partial x}\right|, \ \ \left|\frac{\partial K(x,y,t)}{\partial y}\right| \leq\frac{Ct}{(x^2+y^2)^{2}}.
\end{equation}

\section{Dislocation velocity given by the threshold dynamics method}\label{sec:Ana}

In this section, we examine the dislocation velocity given by the threshold dislocation dynamics method presented in Sec.~\ref{sec:Thresh}, by comparing it with that in the available discrete dislocation dynamics methods~\cite{kubin1992,Ghoniem2000,xiang2003level,Arsenlis2007}.

 In discrete dislocation dynamics~\cite{kubin1992,Ghoniem2000,xiang2003level,Arsenlis2007}, the dislocation velocity is determined from the force on dislocations by
\begin{flalign}
\mathbf v=&M\mathbf f,\\
\mathbf f=&(\pmb \sigma\cdot \mathbf b)\times \pmb \tau,
\end{flalign}
where $M$ is the mobility, $\mathbf f$ is the Peach-Koehler force on a dislocation~\cite{anderson2017theory}, $\pmb \sigma$ is the stress tensor, $\mathbf b$ is the Burgers vector, and $\pmb \tau$ is the unit tangent vector of the dislocation. When the dislocation is in the $xy$ plane and the Burgers vector $\mathbf b=(b,0)$, without the applied stress, the velocity $\mathbf v$ and Peach-Koehler force $\mathbf f$ are in the normal direction of the dislocation with values $v=M\sigma_{13}b$ and $f=\sigma_{13}b$, respectively.
 In this case, the dislocation velocity in the normal direction, due to the stress generated by the dislocations, is~\cite{anderson2017theory,xiang2003level,cai2006non}
 {
\begin{flalign}
v(x,y)=&-\int_\Gamma d\bar{y}\iint_{\mathbf R^2}\frac{x-\xi}{4\pi(1-\nu)[(x-\xi)^2+(y-\eta)^2]^{\frac{3}{2}}}\delta_\varepsilon(\xi-\bar{x},\eta-\bar{y})d\xi d\eta\nonumber\\
&+\int_\Gamma d\bar{x}\iint_{\mathbf R^2}\frac{y-\eta}{4\pi[(x-\xi)^2+(y-\eta)^2]^{\frac{3}{2}}}\delta_\varepsilon(\xi-\bar{x},\eta-\bar{y})d\xi d\eta,\label{eqn:dis-v}
\end{flalign}
}
where $\delta_\varepsilon$ is some two-dimensional regularized delta-function with regularization width $\varepsilon\ll 1$ that represents the dislocation core width. Note that Eq.~\eqref{eqn:dis-v} is in dimensionless form, in which the length has been {stretched} by the unit length of the simulation domain $l_0$ and time by $l_0^2/M\mu b^2$.
Note that the stress and velocity formulation directly given by the dislocation theory~\cite{anderson2017theory}, i.e., without convolution with the regularized delta function $\delta_\varepsilon$ in the above equation, is singular,
 and one of the treatments in dislocation dynamics methods is to use a regularized Dirac delta function, which represents the dislocation core effect, to smooth the integral \cite{xiang2003level,cai2006non}.

We will examine the dislocation velocity given by the threshold dislocation dynamics method presented in Sec.~\ref{sec:Thresh} by comparing it with that in the available discrete dislocation dynamics methods given in Eq.~\eqref{eqn:dis-v}. In fact, using the property $\delta_\varepsilon(-x,-y)=\delta_\varepsilon(x,y)$, the dislocation velocity in Eq.~\eqref{eqn:dis-v} can be written as
\begin{flalign}
v(x,y)=&{ \int_\Gamma-\frac{x-\bar{x}}{4\pi(1-\nu)[(x-\bar{x})^2+(y-\bar{y})^2]^{\frac{3}{2}}}*\delta_\varepsilon d\bar{y}+\frac{y-\bar{y}}{4\pi[(x-\bar{x})^2+(y-\bar{y})^2]^{\frac{3}{2}}}*\delta_\varepsilon d\bar{x}} \nonumber\\
=&-\int_\Gamma d\bar{y}\iint_{\mathbf R^2}\frac{x-\xi-\bar{x}}{4\pi(1-\nu)[(x-\xi-\bar{x})^2+(y-\eta-\bar{y})^2]^{\frac{3}{2}}}\delta_\varepsilon(\xi,\eta)d\xi d\eta\nonumber\\
&+\int_\Gamma d\bar{x}\iint_{\mathbf R^2}\frac{y-\eta-\bar{y}}{4\pi[(x-\xi-\bar{x})^2+(y-\eta-\bar{y})^2]^{\frac{3}{2}}}\delta_\varepsilon(\xi,\eta)d\xi d\eta.\label{eqn:dis-v1}
\end{flalign}
{ Here $``*"$ is the convolution operator in two dimensions with respect to the variable $(\bar{x},\bar{y})$.} In particular, from time $t_n$ to $t_{n+1}=t_n+\Delta t$, the evolution equation  in the threshold dislocation dynamics method is  \eqref{eqn:formulation}, from which we
will analyze the average dislocation velocity given by this formulation. { As explained in the previous section,  $\Delta t$ is the dislocation core parameter in the dimensionless form presented at the end of Sec.~\ref{sec:PN}. }

As already mentioned, unlike those available threshold dynamics methods reviewed in the introduction section, all of which focus on the leading order velocity of the moving front that is proportional to its local curvature on the order of $\log\varepsilon$, here
we will show that the threshold dislocation dynamics method is able to generate the correct nearly singular and long-range dislocation velocity, i.e., in the correct  leading order $\log\varepsilon$ and the next order $O(1)$ contributions to the dislocation velocity, as $\varepsilon\to0$.

Note that in the proposed threshold dislocation dynamics method,
the slow decaying dislocation core function,
which is the dislocation stress kernel  $K(x,y,t)$ in Eq.~\eqref{eqn:Kt} ($1/r^3$-decay as given in Eq.~\eqref{eqn:asymp-K}, where $r$ is the distance to the point on the dislocation),
is unlike the dislocation cores of cut-off \cite{GB1976}, compact support  \cite{xiang2003level,Zhao2012}, or $1/r^7$-decay  \cite{cai2006non} profiles in the available discrete dislocation dynamics methods,
and leads to more complex treatments in the analysis of the leading orders of the resulting dislocation velocity.

\subsection{Velocity  due to stress generated by the dislocations}

We first consider the dislocation velocity  given by the threshold dislocation dynamics method  without applied stress. Recall that in the formulation for the stress generated by the dislocation in Eqs.~\eqref{eqn:sigma13} and \eqref{eqn:sigma13-hat}, the Burgers vectors of the dislocation is $\mathbf b=(b,0,0)$.
In this case, the solution in the evolution step of the method in Eq.~\eqref{eqn:formulation} becomes
\begin{equation}\label{eqn:formulation1}
{u}(x,y,t_{n+1})= K_{\Delta t}*1_{S_n},
\end{equation}
where $K_{\Delta t}$ is given by Eq.~\eqref{eqn:K}.

Without loss of generality, assume $t_n=0$. In addition to ${u}(x,y,\Delta t)= K_{\Delta t}*1_{S_0}$ given by  Eq.~\eqref{eqn:formulation1}, suppose that for $t\in[0, \Delta t]$, $u$ satisfies
\begin{equation}\label{eqn:formulation2}
u(x,y,t)= K*1_{S_0},
\end{equation}
where $K=K(x,y,t)$ is the kernel function given in Eq.~\eqref{eqn:Kt}. Here $t$ serves as the instant dislocation core parameter.
Using the identity $\widehat{f*g}=(2\pi)^2\hat{f}\hat{g}$ and Eq.~\eqref{eqn:formulation2}, we have
$\hat u=\widehat{K*1_{S_0}}=(2\pi)^2\hat K \, \widehat{1}_{S_0}$, and
\begin{equation}
\hat{u}_t=-\frac{1}{2} \left(\frac{k_1^2}{(1-\nu)\|\mathbf{k}\|} +\frac{k_2^2}{\|\mathbf{k}\|}\right)(2\pi)^2\hat K \, \widehat{1}_{S_0}.
\end{equation}
Further using the fact that the inverse Fourier transform of $F(k_1,k_2)=-\frac{1}{2} \left(\frac{k_1^2}{(1-\nu)\|\mathbf{k}\|} +\frac{k_2^2}{\|\mathbf{k}\|}\right)$ is
$g(x,y)=-\frac{\pi}{1-\nu}\frac{\partial }{\partial x}\left(\frac{x}{(x^2+y^2)^\frac{3}{2}}\right)
-\pi\frac{\partial }{\partial y}\left(\frac{y}{(x^2+y^2)^\frac{3}{2}}\right)$,
it can be calculated that
\begin{flalign}
u_t(x,y,t)=&{ \int_\Gamma\frac{x-\bar{x}}{4\pi(1-\nu)[(x-\bar{x})^2+(y-\bar{y})^2]^{\frac{3}{2}}}*K(\bar{x},\bar{y},t)d\bar{y}}\nonumber\\
&{ -\frac{y-\bar{y}}{4\pi[(x-\bar{x})^2+(y-\bar{y})^2]^{\frac{3}{2}}}*K(\bar{x},\bar{y},t)d\bar{x}}\nonumber\\
=&\int_\Gamma d\bar{y}\iint_{\mathbf R^2}\frac{x-\xi-\bar{x}}{4\pi(1-\nu)[(x-\xi-\bar{x})^2+(y-\eta-\bar{y})^2]^{\frac{3}{2}}}K(\xi,\eta,t)d\xi d\eta\nonumber\\
&-\int_\Gamma d\bar{x}\iint_{\mathbf R^2}\frac{y-\eta-\bar{y}}{4\pi[(x-\xi-\bar{x})^2+(y-\eta-\bar{y})^2]^{\frac{3}{2}}}K(\xi,\eta,t)d\xi d\eta.\label{eqn:ut}
\end{flalign}

In order to obtain the asymptotic behavior as $\Delta t\to 0$ of the velocity of the dislocation whose dynamics is given implicitly by ${u}(x,y,\Delta t)= K_{\Delta t}*1_{S_0}$,
we first obtain the asymptotic behavior of $u_t(x,y,t)$ given by Eq.~\eqref{eqn:ut}. We assume that the length of the dislocation $\Gamma$ is $O(1)$.
%
%
%

We write $u_t$ in Eq.~\eqref{eqn:ut} as
\begin{flalign}
u_t(x,y,t)=&\iint_{\mathbf R^2}I(\xi,\eta,x,y,t) K(\xi,\eta,t)d\xi d\eta,
\end{flalign}
where
{
\begin{flalign}
I=&\int_\Gamma \frac{x+\xi-\bar{x}}{4\pi(1-\nu)[(x+\xi-\bar{x})^2+(y+\eta-\bar{y})^2]^{\frac{3}{2}}}d\bar{y}
-\frac{y+\eta-\bar{y}}{4\pi[(x+\xi-\bar{x})^2+(y+\eta-\bar{y})^2]^{\frac{3}{2}}}d\bar{x}.\label{eqn:I000}
\end{flalign}}
Note that here we have changed integration variables from $(\xi,\eta)$ to $(-\xi,-\eta)$ and still use the notations $(\xi,\eta)$, and the property that $K(\xi,\eta,t)$ is an even function with respect to $\xi$ and $\eta$.

We will use the following asymptotic behavior for the integral $I$ in Eq.~\eqref{eqn:I000}, which is the stress (or force, up to a constant factor) generated by the dislocation $\Gamma$  at a point $P$ with distance $d_P$ to the dislocation \cite{GB1976,Zhao2012}:
\begin{flalign}
I=-\left(1+ \frac{\nu\sin^2\alpha}{1-\nu}\right)\frac{{\rm sgn}(P)}{2\pi d_P}+\frac{1+\nu(1-3\sin^2\alpha)}{4\pi(1-\nu)}\kappa\log d_P+O(1), \ \ d_P\rightarrow 0,\label{eqn:I00-asymp}
\end{flalign}
where $\alpha$ is the angle between the Burgers vector and the dislocation line direction at the nearest point to the evaluation point $P$, ${\rm sgn}(P)=1$ or $-1$ when $P$ is on the positive or negative side of the dislocation, and $\kappa$ is the curvature of the dislocation at the nearest point to $P$. Here the point $P=(x+\xi,y+\eta)$ for $(\xi,\eta)$ varying over $\mathbf R^2$

\begin{figure}[!htb]
	\centering
	\includegraphics[width=4.5in]{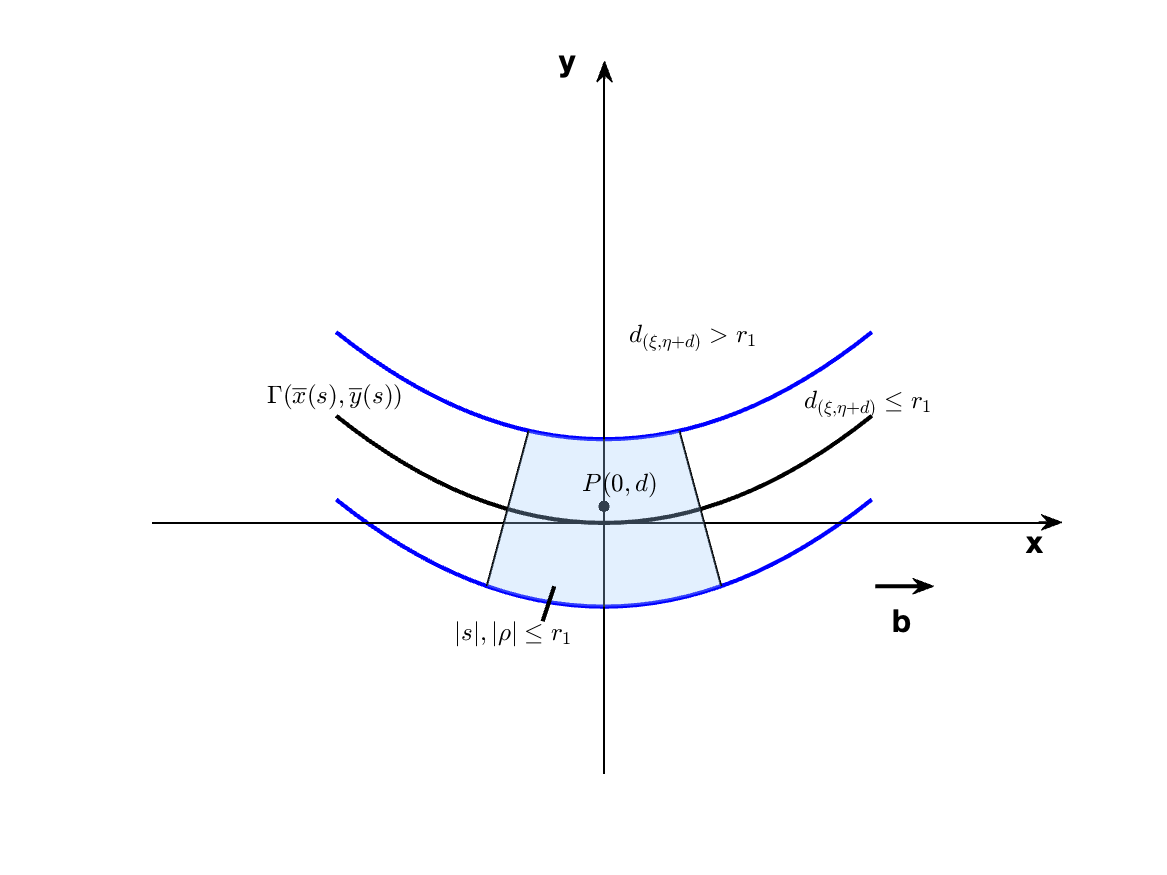}
	\caption{Dislocation $\Gamma$ and different regions for the analysis of dislocation velocity given by the threshold dynamics method. The dislocation core region is the tube between the two blue curves, and the
 black line in between is the dislocation $\Gamma$.  The region
 $|s|, |\rho|\leq r_1$ is the shaded region that contains the origin and the point $P(0,d)$. The set $D_0$ consists of all  $(\xi,\eta)$ such that $(\xi,\eta+d)$ is inside the dislocation core, i.e., $d_{(\xi,\eta+d)}\leq r_1$, while the set $D_1$ consists of all  $(\xi,\eta)$ such that $(\xi,\eta+d)$ is outside the dislocation core, i.e., $d_{(\xi,\eta+d)}>r_1$.}
	\label{fig:curve}
\end{figure}

 Suppose that { the point $(0,0)$ is on the  dislocation $\Gamma$  and the dislocation is in the $+x$ direction at $(0,0)$}, see Fig.~\ref{fig:curve}. That is, the dislocation is screw at the point $(0,0)$.
Consider $u_t$ at a point $(0,d)$, where $d<\Delta t$.

{ In this case, the point $(x,y)=(0,d)$ in the integral $I$ in  Eq.~\eqref{eqn:I000}, and we have}
\begin{flalign}
I(\xi,\eta,d,t)=&\int_\Gamma \frac{\xi-\bar{x}}{4\pi(1-\nu)[(\xi-\bar{x})^2+(d+\eta-\bar{y})^2]^{\frac{3}{2}}}d\bar{y}
-\frac{d+\eta-\bar{y}}{4\pi[(\xi-\bar{x})^2+(d+\eta-\bar{y})^2]^{\frac{3}{2}}}d\bar{x};\label{eqn:I00}
\end{flalign}
{ and in the asymptotic behavior in Eq.~\eqref{eqn:I00-asymp}, as discussed above,
 the point $P=(\xi,\eta+d)$ for $(\xi,\eta)$ varying over $\mathbf R^2$, and $\alpha=0$ for the point $P=(0,d)$.}


Case 1. $D_1=\{(\xi,\eta): d_{(\xi,\eta+d)}>r_1\}$, { where $r_1$ is a small number and $r_1>2d$.
Here following the notation $d_P$ defined above, the notation $d_{(\xi,\eta)}$ is the distance from the point $(\xi,\eta)$ to the dislocation $\Gamma$. }

In this case, we have $d_{(\xi,\eta)}\geq d_{(\xi,\eta+d)}-d>\frac{r_1}{2}$.
In particular, for $(\xi,\eta)\in D_1$, we have $\sqrt{\xi^2+\eta^2}> \frac{r_1}{2}$.
{ We also have that for $(\xi,\eta)\in D_1$, the denominator in the integrand of $I(\xi,\eta,d,t)$ in Eq.~\eqref{eqn:I00}, excluding the constant factor, is $$\|(\xi,\eta+d)-(\bar{x},\bar{y})\|\geq \|(\xi,\eta)-(\bar{x},\bar{y})\|-d\geq r_1-d>\frac{r_1}{2}.$$
Thus, for $(\xi,\eta)\in D_1$, we have
$$|I(\xi,\eta,d,t)|\leq \frac{CL}{r_1^2},$$
where $L$ is the total length of the dislocation.}

Further using the upper bound of $K(\xi,\eta,t)$ in Eq.~\eqref{eqn:asymp-K} and accordingly,
\begin{equation}
\iint_{D_1}|K(\xi,\eta,t)|d\xi d\eta\leq \iint_{\sqrt{\xi^2+\eta^2}>r_1}|K(\xi,\eta,t)|d\xi d\eta\leq \dfrac{Ct}{r_1}, \label{eqn:K-out0}
\end{equation}
 we have
$$\left|\iint_{D_1}I(\xi,\eta,d,t) K(\xi,\eta,t)d\xi d\eta\right|\leq \frac{CL}{r_1^2}\iint_{\sqrt{\xi^2+\eta^2}>r_1}|K(\xi,\eta,t)|d\xi d\eta\leq \dfrac{CLt}{r_1^3}.$$
Thus
\begin{flalign}
\iint_{D_1}I(\xi,\eta,d,t) K(\xi,\eta,t)d\xi d\eta
=O\left(\frac{t}{r_1^3} \right).\label{eqn:outsidecore}
\end{flalign}

Case 2. $D_0=\{(\xi,\eta): d_{(\xi,\eta+d)}\leq r_1\}$. 
Using the asymptotic behavior in Eq.~\eqref{eqn:I00-asymp}, we only need to consider the contributions of $O(1/d_P)$, $O(\log d_P)$, and $O(1)$.

{ Assume that the dislocation is $(\bar{x}(s),\bar{y}(s))$, where $s$ be the arclength parameter of the dislocation $\Gamma$. The dislocation core region $D_0$ can be written as \begin{equation}
(\xi,\eta+d)=(\bar{x}(s),\bar{y}(s))+\rho(-\bar{y}'(s),\bar{x}'(s)), \ \ \ (s,\rho)\in [-\frac{L}{2},\frac{L}{2}]\times[-r_1,r_1]. \label{eqn:core-parametrization}
\end{equation}
{ By  Eq.~\eqref{eqn:core-parametrization}, the closest point on the dislocation to the point $(\xi,\eta+d)$ is
$(\bar{x}(s),\bar{y}(s))$, and the distance between them is $\rho$. This is because the segment between these two points is perpendicular to the tangent vector $(\bar{x}'(s),\bar{y}'(s))$ at the point $(\bar{x}(s),\bar{y}(s))$.}

   Near the point  $(0,0)$, the dislocation $\Gamma$ in the local canonical form \cite{curve-canonical} is
\begin{equation}\label{eqn:parametrization000}
\bar{x}(s)=s+O(s^3), \ \ \ \bar{y}(s)=\frac{\kappa}{2}s^2+O(s^3), \ \ s\rightarrow 0,
\end{equation}
where  $\kappa$ is the curvature of the dislocation at point $(0,0)$.}

(i) \underline{$O(1)$ contribution}.

Denote the $O(1)$ contribution in Eq.~\eqref{eqn:I00-asymp} as $-A(\xi,\eta)$, where $A(\xi,\eta)$ is a bounded function with bounded partial derivatives. { Using the upper bound of $K(\xi,\eta,t)$ in Eq.~\eqref{eqn:asymp-K}, we have $ \iint_{\sqrt{\xi^2+\eta^2}>\frac{r_1}{2}}|K(\xi,\eta,t)|d\xi d\eta\leq \dfrac{Ct}{r_1}$. Thus,}
\begin{flalign}
&\iint_{D_0}A(\xi,\eta) K(\xi,\eta,t)d\xi d\eta \nonumber\\
=&\left(\iint_{D_0: \sqrt{\xi^2+\eta^2}\leq \frac{r_1}{2}}+\iint_{D_0: \sqrt{\xi^2+\eta^2}\geq \frac{r_1}{2}}\right)A(\xi,\eta) K(\xi,\eta,t)d\xi d\eta \nonumber\\
=&\iint_{D_0: \sqrt{\xi^2+\eta^2}\leq \frac{r_1}{2}}A(\xi,\eta) K(\xi,\eta,t)d\xi d\eta +O\left(\frac{t}{r_1}\right)\nonumber\\
=&A(0,0)\iint_{D_0: \sqrt{\xi^2+\eta^2}\leq \frac{r_1}{2}}K(\xi,\eta,t)d\xi d\eta+O(r_1) +O\left(\frac{t}{r_1}\right)\nonumber\\
=&A(0,0)\iint_{\sqrt{\xi^2+\eta^2}\leq \frac{r_1}{2}}K(\xi,\eta,t)d\xi d\eta+O(r_1) +O\left(\frac{t}{r_1}\right)\nonumber\\
=&A(0,0)\left(1-\iint_{\sqrt{\xi^2+\eta^2}\geq \frac{r_1}{2}}K(\xi,\eta,t)d\xi d\eta\right)+O(r_1) +O\left(\frac{t}{r_1}\right)\nonumber\\
=&A(0,0)+O(r_1) +O\left(\frac{t}{r_1}\right).\label{eqn:0-g}
\end{flalign}
Note that here $A(0,0)$ is the $O(1)$ contribution in Eq.~\eqref{eqn:I00-asymp} at the point $P=(0,d)$.

(ii) \underline{$O(\log d_P)$ contribution}.

 Using  the parametrization in Eq.~\eqref{eqn:core-parametrization} { and the local canonical form in Eq.~\eqref{eqn:parametrization000}},
we have $\xi=s-\kappa \rho s+{ O(s^3+\rho s^2)}$, for $s,\rho\in[-r_1,r_1]$, where $\kappa$ is the curvature of the dislocation at $(0,0)$, i.e. $s=0$. { Thus,  for $(\xi,\eta)\in D_0$ with $|s|>r_1$ and small enough $r_1$, we have $\sqrt{\xi^2+\eta^2}\geq O(r_1)$.}

For a point $(\xi,\eta)\in D_0$ with   $s,\rho\in[-r_1,r_1]$,
 { using Eq.~\eqref{eqn:core-parametrization} and the local canonical form in Eq.~\eqref{eqn:parametrization000}},
 it can be calculated that the signed distance ${\rm sgn}(P)d_P$ from the point $(\xi,\eta+d)$ to the dislocation is
\begin{equation}
{\rm sgn}(P)d_P=\rho=\eta+d-\frac{\kappa}{2}\xi^2+{ O(s^3+\rho s^2)}. \label{eqn:d-P}
\end{equation}



{ We write the $O(\log d_P)$ contribution as
\begin{flalign*}
&\iint_{D_0} {\textstyle \frac{1+\nu(1-3\sin^2\alpha(s))}{4\pi(1-\nu)}}\kappa(s)\log d_P \cdot K(\xi,\eta,t)d\xi d\eta \nonumber\\
=&\iint_{D_0: |s|\leq r_1} {\textstyle \frac{1+\nu(1-3\sin^2\alpha(s))}{4\pi(1-\nu)}}\kappa(s)\log d_P \cdot K(\xi,\eta,t)d\xi d\eta \nonumber\\ &+\iint_{D_0: |s|>r_1} {\textstyle \frac{1+\nu(1-3\sin^2\alpha(s))}{4\pi(1-\nu)}}\kappa(s)\log d_P \cdot K(\xi,\eta,t)d\xi d\eta \nonumber\\
\equiv&J_1+J_2.
\end{flalign*}

For the contribution $J_2$, which is for $(\xi,\eta)\in D_0$ with $|s|>r_1$, we have shown that $\sqrt{\xi^2+\eta^2}\geq O(r_1)$, and then using the upper bound of $K(\xi,\eta,t)$ in Eq.~\eqref{eqn:asymp-K}, we have $K(\xi,\eta,t)=O\left(\frac{t}{r_1^3}\right)$. Further using the parametrization in Eq.~\eqref{eqn:core-parametrization}, we have
\begin{flalign*}
J_2=&\iint_{D_0: |s|>r_1} {\textstyle \frac{1+\nu(1-3\sin^2\alpha(s))}{4\pi(1-\nu)}}\kappa(s)\log d_P \cdot K(\xi,\eta,t)d\xi d\eta \nonumber\\
=& \iint_{D_0: |s|>r_1} \log d_P \cdot O\left(\frac{t}{r_1^3}\right) d\xi d\eta \nonumber\\
=&\int_{-r_1}^{r_1} d\rho \, \log |\rho| \cdot(1-\kappa \rho)\int_{r_1<|s|\leq \frac{L}{2}}  O\left(\frac{t}{r_1^3}\right) ds \nonumber\\
=& O\left(\frac{t}{r_1^2}\log r_1\right).\nonumber\\
\end{flalign*}

Now we consider the contribution $J_1$. From Eqs.~\eqref{eqn:core-parametrization} and \eqref{eqn:parametrization000}, we have
$\xi=s-\kappa \rho s+\rho(\kappa s-\bar{y}'(s))+O(s^3)$, $\eta+d=\rho+\frac{\kappa}{2}s^2+\rho(\bar{x}'(s)-1)+O(s^3)$.
Let
\begin{flalign*}
&\xi_0=s-\kappa \rho s+\rho(\kappa s-\bar{y}'(s))=s-\kappa \rho s+O(\rho s^2),\\
&\eta_0+d=\rho+\frac{\kappa}{2}s^2+\rho(\bar{x}'(s)-1)=\rho+\frac{\kappa}{2}s^2+O(\rho s^2).
\end{flalign*}
We have
$\xi=\xi_0+O(s^3)=\xi_0+O(\xi_0^3)$, $\eta=\eta_0+O(s^3)=\eta_0+O(\xi_0^3)$, and
\begin{flalign*}
\rho=\eta_0+d-\frac{\kappa}{2}s^2+O(\rho s^2)=\eta_0+d-\frac{\kappa}{2}\xi_0^2+O(\rho s^2).
\end{flalign*}
This gives $\rho=\left(\eta_0+d-\frac{\kappa}{2}\xi_0^2\right)(1+O(r_1^2))$. Thus, we have
\begin{flalign}
J_1=&\iint_{D_0: |s|\leq r_1} {\textstyle \frac{1+\nu(1-3\sin^2\alpha(s))}{4\pi(1-\nu)}}\kappa(s)\log d_P \cdot K(\xi,\eta,t)d\xi d\eta \nonumber\\
=&\iint_{\scriptsize \begin{array}{l}\xi_0=s-\kappa \rho s+\rho(\kappa s-\bar{y}'(s)), \\ \eta_0+d=\rho+\frac{\kappa}{2}s^2+\rho(\bar{x}'(s)-1), \\ |s|, |\rho|\leq r_1\end{array}}
 {\textstyle \frac{1+\nu}{4\pi(1-\nu)}}\kappa
 \log \left(\left|\eta_0+d-\frac{\kappa}{2}\xi_0^2\right|\Big(1+O(r_1^2)\Big) \right)\nonumber\\
 &\cdot K\Big(\xi_0+O(\xi_0^3),\eta_0+O(\xi_0^3),t\Big) \cdot (1+O(r_1))d\xi_0 d\eta_0\nonumber\\
=&\iint_{\scriptsize \begin{array}{l}\xi_0=s-\kappa \rho s+\rho(\kappa s-\bar{y}'(s)), \\ \eta_0+d=\rho+\frac{\kappa}{2}s^2+\rho(\bar{x}'(s)-1), \\ |s|, |\rho|\leq r_1\end{array}}
 {\textstyle \frac{1+\nu}{4\pi(1-\nu)}}\kappa
 \log \left|\eta_0+d-\frac{\kappa}{2}\xi_0^2\right|
 \cdot K\Big(\xi_0,\eta_0,t\Big) \cdot (1+O(r_1))d\xi_0 d\eta_0\nonumber\\&+O(r_1^2)+O(t r_1\log r_1)\nonumber\\
=&{\textstyle \frac{1+\nu}{4\pi(1-\nu)}}\kappa
\iint_{\mathbf R^2}
  \log \left|\eta_0+d-\frac{\kappa}{2}\xi_0^2\right|
 \cdot K\Big(\xi_0,\eta_0,t\Big) \cdot (1+O(r_1))d\xi_0 d\eta_0\nonumber\\
 &+O(r_1^2)+O(t r_1\log r_1)+O\left(\frac{t}{r_1^2}\log r_1\right)+O\left(\frac{t}{r_1}\log r_1\right).\nonumber
\end{flalign}
Here in the second equation,
the factor $1+O(r_1)$ at the end of the integrand comes from the change of variables from $(\xi,\eta)$ to $(\xi_0,\eta_0)$,
and the approximation of the angle $\alpha$ dependent prefactor by its value at the origin.
In the third equation, we have used the bound of partial derivatives of $K$ in Eq.~\eqref{eqn:asymp-K1} to obtain the error $O(t r_1\log r_1)$ when $K\Big(\xi_0+O(\xi_0^3),\eta_0+O(\xi_0^3),t\Big)$ is replaced by $K\Big(\xi_0,\eta_0,t\Big)$, and the error $O(r_1^2)$ comes from the $O(r_1^2)$ term inside the logarithm.
In the fourth equation, we extend the integration domain to the entire $\mathbf R^2$, in which the region $D_0: |s|>r_1$ added in this step gives the error of  $O\left(\frac{t}{r_1^2}\log r_1\right)$ as the calculation of $J_2$, and the region $D_1$ added in this step gives the error of $O\left(\frac{t}{r_1}\log r_1\right)$ using the fact that $\sqrt{\xi^2+\eta^2}\geq \frac{r_1}{2}$ proved in Case 1 and accordingly $\sqrt{\xi_0^2+\eta_0^2}\geq \frac{r_1}{4}$ for small enough $r_1$ together with Eq.~\eqref{eqn:asymp-K}.

Combining the results of $J_1$ and $J_2$, we have
\begin{flalign}
&\iint_{D_0} {\textstyle \frac{1+\nu(1-3\sin^2\alpha(s))}{4\pi(1-\nu)}}\kappa(s)\log d_P \cdot K(\xi,\eta,t)d\xi d\eta \nonumber\\
=&{\textstyle \frac{1+\nu}{4\pi(1-\nu)}}\kappa
\iint_{\mathbf R^2}
  \log \left|\eta+d-\frac{\kappa}{2}\xi^2\right|
 \cdot K(\xi,\eta,t) \cdot (1+O(r_1))d\xi d\eta+O(r_1^2)+O\left(\frac{t}{r_1^2}\log r_1\right).
\label{eqn:log-term}
\end{flalign}
Here we still use variables $(\xi,\eta)$ instead of $(\xi_0,\eta_0)$  for simplicity of notations, and combine all the errors.

}

%

%
%

The result in Eq.~\eqref{eqn:log-term} can be further simplified. Using $K(\xi,\eta,t)d\xi d\eta=K(\frac{\xi}{t},\frac{\eta}{t},1)d\frac{\xi}{t} d\frac{\eta}{t}$,
and denoting $\xi_1=\frac{\xi}{t}$, $\eta_1=\frac{\eta}{t}$,
 the leading order term in Eq.~\eqref{eqn:log-term} is
\begin{flalign}
&\iint_{\mathbf R^2} \log \left|\eta+d-\frac{\kappa}{2}\xi^2\right| \cdot K(\xi,\eta,t)d\xi d\eta\nonumber\\
=&\iint_{\mathbf R^2} \log \left|t\eta_1+d-\frac{\kappa}{2}t^2\xi_1^2\right| \cdot K(\xi_1,\eta_1,1)d\xi_1 d\eta_1\nonumber\\
=&\log t+\iint_{\mathbf R^2} \log \left|\eta_1+\frac{d}{t}-\frac{\kappa}{2}t\xi_1^2\right| \cdot K(\xi_1,\eta_1,1)d\xi_1 d\eta_1\nonumber\\
=&\log t+\iint_{\sqrt{\xi_1^2+\eta_1^2}\leq t^{-\frac{1}{3}}} \log \left|\eta_1+\frac{d}{t}-\frac{\kappa}{2}t\xi_1^2\right| \cdot K(\xi_1,\eta_1,1)d\xi_1 d\eta_1\nonumber\\
&+\iint_{\sqrt{\xi_1^2+\eta_1^2}>t^{-\frac{1}{3}}, \ |\eta_1+\frac{d}{t}-\frac{\kappa}{2}t\xi_1^2|> \frac{1}{2}t^{-\frac{1}{3}}} \log \left|\eta_1+\frac{d}{t}-\frac{\kappa}{2}t\xi_1^2\right| \cdot K(\xi_1,\eta_1,1)d\xi_1 d\eta_1\nonumber\\
&+\iint_{\sqrt{\xi_1^2+\eta_1^2}>t^{-\frac{1}{3}},\ |\eta_1+\frac{d}{t}-\frac{\kappa}{2}t\xi_1^2|\leq \frac{1}{2}t^{-\frac{1}{3}}} \log \left|\eta_1+\frac{d}{t}-\frac{\kappa}{2}t\xi_1^2\right| \cdot K(\xi_1,\eta_1,1)d\xi_1 d\eta_1\nonumber\\
=&\log t+\iint_{\sqrt{\xi_1^2+\eta_1^2}\leq t^{-\frac{1}{3}}} \log \left|\eta_1+\frac{d}{t}-\frac{\kappa}{2}t\xi_1^2\right| \cdot K(\xi_1,\eta_1,1)d\xi_1 d\eta_1+O\left(t^{\frac{1}{3}}\log t\right)\nonumber\\
=&\log t+\iint_{\sqrt{\xi_2^2+\eta_2^2}\leq t^{-\frac{1}{3}}} \log \left|\eta_2+\frac{d}{t}\right| \cdot K\left(\xi_2,\eta_2+\frac{\kappa}{2}t\xi_2^2,1\right)d\xi_2 d\eta_2+O\left(t^{\frac{1}{3}}\log t\right)\nonumber\\
=&\log t+\iint_{\sqrt{\xi_2^2+\eta_2^2}\leq t^{-\frac{1}{3}}} \log \left|\eta_2+\frac{d}{t}\right| \cdot K\left(\xi_2,\eta_2,1\right)d\xi_2 d\eta_2+O\left(t^{\frac{1}{3}}\log t\right)\nonumber\\
&+\iint_{\sqrt{\xi_2^2+\eta_2^2}\leq  1} \log \left|\eta_2+\frac{d}{t}\right| \cdot \left[K\left(\xi_2,\eta_2+\frac{\kappa}{2}t\xi_2^2,1\right)-K(\xi_2,\eta_2,1)\right]d\xi_2 d\eta_2\nonumber\\
&+\iint_{1\leq\sqrt{ \xi_2^2+\eta_2^2}\leq t^{-\frac{1}{3}}} \log \left|\eta_2+\frac{d}{t}\right| \cdot \left[K\left(\xi_2,\eta_2+\frac{\kappa}{2}t\xi_2^2,1\right)-K(\xi_2,\eta_2,1)\right]d\xi_2 d\eta_2\nonumber\\
=&\log t+\iint_{\sqrt{\xi_2^2+\eta_2^2}\leq t^{-\frac{1}{3}}} \log \left|\eta_2+\frac{d}{t}\right| \cdot K\left(\xi_2,\eta_2,1\right)d\xi_2 d\eta_2+O\left(t^{\frac{1}{3}}\log t\right)\nonumber\\
&+O\left(t^{\frac{1}{3}}\right)+O\left(t^{\frac{1}{3}}\right)\left(1+t^{\frac{2}{3}}\log t\right)\nonumber\\
=&\log t+\iint_{\mathbf R^2} \log \left|\eta_2+\frac{d}{t}\right| \cdot K\left(\xi_2,\eta_2,1\right)d\xi_2 d\eta_2+O\left(t^{\frac{1}{3}}\log t\right). \label{eqn:log-term1}
\end{flalign}
Here the first error term $O\left(t^{\frac{1}{3}}\log t\right)$ is obtained by using the bound of $K$ in Eq.~\eqref{eqn:asymp-K}.
The change of integration variables is from $(\xi_1,\eta_1)$ to $\xi_2=\xi_1$ and $\eta_2=\eta_1-\frac{\kappa}{2}t\xi_1^2$, and we have $t\xi^2=O(t^\frac{1}{3})$ when $\sqrt{\xi_1^2+\eta_1^2}\leq t^{-\frac{1}{3}}$.
The integrals containing the difference $K\left(\xi_2,\eta_2+\frac{\kappa}{2}t\xi_2^2,1\right)-K(\xi_2,\eta_2,1)$ are estimated by using the mean value theorem and the bound of $\frac{\partial K}{\partial y}$ in Eq~\eqref{eqn:asymp-K1}. 

Summarizing Eqs.~\eqref{eqn:log-term} and \eqref{eqn:log-term1}, and using
\begin{flalign}
&\iint_{\mathbf R^2} \log \left|\eta_2+\frac{d}{t}\right| K(\xi_2,\eta_2,1)d\xi_2 d\eta_2\nonumber\\
=&\frac{1}{(2\pi)^2}\iint_{\mathbf R^4}  \log \left|\eta_2+\frac{d}{t}\right|  e^{-\frac{1}{2} \left(\frac{k_1^2}{(1-\nu)\|\mathbf{k}\|} +\frac{k_2^2}{\|\mathbf{k}\|}\right)}e^{i(k_1\xi_2+k_2\eta_2)} \mathrm{d}k_1\mathrm{d}k_2 d\xi_2 d\eta_2\nonumber\\
=&\frac{1}{2\pi}\iint_{\mathbf R^2}  \log \left|\eta_2+\frac{d}{t}\right|  e^{-\frac{1}{2} |k_2|}e^{ik_2\eta_2} d k_2 d\eta_2\nonumber\\
=&\frac{1}{2\pi}\int_{-\infty}^{\infty}    \frac{\log \left|\eta_2+\frac{d}{t}\right|}{(1/2)^2+\eta_2^2}d\eta_2\nonumber\\
=&\log C_2
+\log\sqrt{\frac{4d^2}{t^2}+1},
\end{flalign}
where
\begin{equation}
\log C_2\equiv \frac{1}{2\pi}\int_{-\infty}^{\infty}    \frac{\log \left|\eta_2\right|}{(1/2)^2+\eta_2^2}d\eta_2{,}
\end{equation}
we have
\begin{flalign}
&\iint_{D_0} \frac{1+\nu(1-3\sin^2\alpha(s))}{4\pi(1-\nu)}\kappa(s)\log d_P \cdot K(\xi,\eta,t)d\xi d\eta \nonumber\\
=& \frac{1+\nu}{4\pi(1-\nu)}\kappa \log (C_2\sqrt{4d^2+t^2} )+O\left(\frac{t}{r_1^2}\log r_1\right) +O(r_1\log t)
+O\left(t^{\frac{1}{3}}\log t\right).
\end{flalign}
This formulation holds when the point $(0,0)$ on the dislocation is screw.

In a general case, where the angle between the line direction of the dislocation and the Burgers vector is $\alpha$ at the point $(0,0)$ being considered, { for the $\log d_P$ term in  Eq.~\eqref{eqn:I00-asymp}, similar calculations give that}
\begin{flalign}
&\iint_{D_0} \frac{1+\nu(1-3\sin^2\alpha(s))}{4\pi(1-\nu)}\kappa(s)\log d_P \cdot K(\xi,\eta,t)d\xi d\eta \nonumber\\
=& \frac{1+\nu(1-3\sin^2\alpha)}{4\pi(1-\nu)}\kappa \log\left[
{\textstyle \left(1+\frac{\nu}{1-\nu}\sin^2\alpha\right)} C_2\sqrt{4d^2+t^2}\right]\nonumber\\
&+O\left(\frac{t}{r_1^2}\log r_1\right) +O(r_1\log t)
+O\left(t^{\frac{1}{3}}\log t\right).\label{eqn:log-d-g}
\end{flalign}

(iii) \underline{$O(1/d_P)$ contribution}.

Using the formula of the distance to the dislocation $d_P$ in Eq.~\eqref{eqn:d-P}, similar to the calculation of the $O(\log d_P)$ contribution given above,  we have
\begin{flalign}
&\iint_{D_0} \left(1+{\textstyle \frac{\nu\sin^2\alpha(s)}{1-\nu}}\right)\frac{{\rm sgn}(P)}{2\pi d_P} K(\xi,\eta,t)d\xi d\eta \nonumber\\
=&\iint_{D_0: |s|\leq r_1}\left(1+{\textstyle \frac{\nu\sin^2\alpha(s)}{1-\nu}}\right)\frac{1}{2\pi \rho} K(\xi,\eta,t)d\xi d\eta \nonumber\\
& +\iint_{D_0: |s|>r_1} \left(1+{\textstyle \frac{\nu\sin^2\alpha(s)}{1-\nu}}\right)\frac{1}{2\pi \rho} K(\xi,\eta,t)d\xi d\eta \nonumber\\
=&\frac{1}{2\pi}\iint_{|s|, |\rho|\leq r_1}
 \frac{1}{\rho}\cdot K(\xi,\eta,t)d\xi d\eta\nonumber\\
& + \frac{1}{2\pi}\int_{r_1<|s|\leq \frac{L}{2}} ds \left(1+{\textstyle \frac{\nu\sin^2\alpha(s)}{1-\nu}}\right)\int_{-r_1}^{r_1} \frac{1}{\rho} K(\xi,\eta,t)(1-\kappa(s)\rho)d\rho   \nonumber\\
=&\frac{1}{2\pi}\iint_{|s|, |\rho|\leq r_1}
 \frac{1}{\rho}\cdot K(\xi,\eta,t)d\xi d\eta
+ \frac{1}{2\pi}\int_{r_1<|s|\leq \frac{L}{2}} ds\left(1+{\textstyle \frac{\nu\sin^2\alpha(s))}{1-\nu}}\right)\nonumber\\
&\ \ \ \ \ \cdot \left[\int_{0}^{r_1}{\textstyle \frac{ K(\xi(s,\rho),\eta(s,\rho),t)-K(\xi(s,-\rho),\eta(s,-\rho),t)}{\rho}}d\rho -\int_{-r_1}^{r_1}  K(\xi,\eta,t)\kappa(s) d\rho\right]\nonumber\\
=& \frac{1}{2\pi}\iint_{|s|, |\rho|\leq r_1}
 \frac{1}{\rho}\cdot K(\xi,\eta,t)d\xi d\eta +O\left(\frac{t}{r_1^2}\right)\nonumber\\
=& \frac{1}{2\pi}\iint_{\mathbf R^2}
 \frac{1}{\rho}\cdot K(\xi,\eta,t)d\xi d\eta +O\left(\frac{t}{r_1^2}\right).\label{eqn:1/d-term-c}
\end{flalign}

Using $K(\xi,\eta,t)d\xi d\eta=K(\frac{\xi}{t},\frac{\eta}{t},1)d\frac{\xi}{t} d\frac{\eta}{t}$,
and denoting $\xi_1=\frac{\xi}{t}$, $\eta_1=\frac{\eta}{t}$, we have
\begin{flalign}
&\frac{1}{2\pi}\iint_{\mathbf R^2}
 \frac{1}{\rho}\cdot K(\xi,\eta,t)d\xi d\eta \nonumber\\
=&\frac{1}{2\pi}\iint_{\sqrt{\xi^2+\eta^2}\leq t^{\frac{2}{5}}}
 \frac{1}{\rho}\cdot K(\xi,\eta,t)d\xi d\eta +O( t^{\frac{1}{5}})\nonumber\\
=&\frac{1}{2\pi}\iint_{\sqrt{\xi^2+\eta^2}\leq t^{\frac{2}{5}}}
 \frac{1}{\eta+d-\frac{\kappa}{2}\xi^2+O(\xi^3+\eta^3)}\cdot K(\xi,\eta,t)d\xi d\eta +O( t^{\frac{1}{5}})\nonumber\\
=&\frac{1}{2\pi t}\iint_{\sqrt{\xi_1^2+\eta_1^2}\leq t^{-\frac{3}{5}}}
 \frac{1}{\eta_1+\frac{d}{t}-\frac{\kappa}{2}t\xi_1^2+O(t^2\xi_1^3+t^2\eta_1^3)}\cdot K(\xi_1,\eta_1,1)d\xi_1 d\eta_1 +O( t^{\frac{1}{5}}) \nonumber\\
=&\frac{1}{2\pi t}\iint_{\sqrt{\xi_2^2+\eta_2^2}\leq t^{-\frac{3}{5}}}
 \frac{1}{\eta_2+\frac{d}{t}-\frac{\kappa}{2}t\xi_2^2}\cdot K(\xi_2,\eta_2+O(t^2\xi_2^3+t^2\eta_2^3),1)d\xi_2 d\eta_2+O(t^{\frac{1}{5}}) \nonumber\\
=&\frac{1}{2\pi t}\iint_{\sqrt{\xi_2^2+\eta_2^2}\leq t^{-\frac{3}{5}}}
 \frac{1}{\eta_2+\frac{d}{t}-\frac{\kappa}{2}t\xi_2^2}\cdot K(\xi_2,\eta_2,1)d\xi_2 d\eta_2+O(t)+O(t^{\frac{1}{5}}) \nonumber\\
=&\frac{1}{2\pi t}\iint_{\sqrt{\xi_2^2+\eta_2^2}\leq t^{-\frac{3}{5}}}
 \frac{1}{\eta_2+\frac{d}{t}}\cdot K(\xi_2,\eta_2,1)d\xi_2 d\eta_2+C_1^{\kappa,d/t,t}+O(t^{\frac{1}{5}}) \nonumber\\
=&\frac{1}{2\pi t}\iint_{\mathbf R^2}
 \frac{1}{\eta_2+\frac{d}{t}}\cdot K(\xi_2,\eta_2,1)d\xi_2 d\eta_2+C_1^{\kappa,d/t,t}+O(t^{\frac{1}{5}}), \label{eqn:1/d-term-c1}
\end{flalign}
where
\begin{flalign}
C_1^{\kappa,\alpha,d/t,t}=&\frac{1}{2\pi t}\left(1+{\textstyle \frac{\nu\sin^2\alpha}{1-\nu}}\right)\iint_{\sqrt{\xi_2^2+\eta_2^2}\leq t^{-\frac{3}{5}}}
 \left(\frac{1}{\eta_2+\frac{d}{t}-\frac{\kappa}{2}t\xi_2^2}- \frac{1}{\eta_2+\frac{d}{t}}\right)\cdot K_\alpha(\xi_2,\eta_2,1)d\xi_2 d\eta_2\nonumber\\
=&\frac{1}{2\pi}\left(1+{\textstyle \frac{\nu\sin^2\alpha}{1-\nu}}\right)\iint_{\sqrt{\xi_2^2+\eta_2^2}\leq t^{-\frac{3}{5}}}
 \frac{\frac{\kappa}{2}\xi_2^2}{\left(\eta_2+\frac{d}{t}-\frac{\kappa}{2}t\xi_2^2\right)\left(\eta_2+\frac{d}{t} \right)}\cdot K_\alpha(\xi_2,\eta_2,1)d\xi_2 d\eta_2\nonumber\\
 =&O(1),
\end{flalign}
with $K_\alpha(x,y,t)=K(x\cos\alpha-y\sin\alpha,x\sin\alpha+y\cos\alpha,t)$ and the angle $\nu=0$ in Eq.~\eqref{eqn:1/d-term-c1}.
Here in the calculation of Eq.~\eqref{eqn:1/d-term-c1}, we have used the upper bounds of $\frac{\partial K}{\partial x}$ and $\frac{\partial K}{\partial y}$ in Eq.~\eqref{eqn:asymp-K1}.  The change of integration variables is from $(\xi_1,\eta_1)$ to $\xi_2=\xi_1$ and $\eta_2=\eta_1+O(t^2\xi_1^3+t^2\eta_1^3)$, and we have $t^2\xi_1^3+t^2\eta_1^3=O(t^\frac{1}{5})$ when $\sqrt{\xi_1^2+\eta_1^2}\leq t^{-\frac{3}{5}}$.


Summarizing Eqs.~\eqref{eqn:1/d-term-c} and \eqref{eqn:1/d-term-c1}, and using
\begin{flalign}
&\iint_{\mathbf R^2} \frac{1}{\eta_2+\frac{d}{t}} K(\xi_2,\eta_2,1)d\xi_2 d\eta_2\nonumber\\
=&\frac{1}{(2\pi)^2}\iint_{\mathbf R^4}  \frac{1}{\eta_2+\frac{d}{t}}   e^{-\frac{1}{2} \left(\frac{k_1^2}{(1-\nu)\|\mathbf{k}\|} +\frac{k_2^2}{\|\mathbf{k}\|}\right)}e^{i(k_1\xi_2+k_2\eta_2)} \mathrm{d}k_1\mathrm{d}k_2 d\xi_2 d\eta_2\nonumber\\
=&\frac{1}{2\pi}\iint_{\mathbf R^2}  \frac{1}{\eta_2+\frac{d}{t}}   e^{-\frac{1}{2} |k_2|}e^{ik_2\eta_2} d k_2 d\eta_2\nonumber\\
=&\frac{1}{2\pi}\int_{-\infty}^{\infty}    \frac{1}{\eta_2+\frac{d}{t}} \frac{1}{\eta_2^2+(1/2)^2}d\eta_2\nonumber\\
=&\frac{td}{\left(\frac{t}{2}\right)^2+d^2},
\end{flalign}
we have
\begin{flalign}
\iint_{D_0} \left(1+{\textstyle \frac{\nu\sin^2\alpha(s)}{1-\nu}}\right)\frac{{\rm sgn}(P)}{2\pi d_P} K(\xi,\eta,t)d\xi d\eta
=\frac{1}{2\pi }\frac{d}{\left(\frac{t}{2}\right)^2+d^2}+C_1^{\kappa,\alpha,d/t,t}+O(t^{\frac{1}{5}})+O\left(\frac{t}{r_1^2}\right).
\end{flalign}
This formulation holds when the point $(0,0)$ on the dislocation is screw.

In a general case, where the angle between the line direction of the dislocation and the Burgers vector is $\alpha$ at the point $(0,0)$ being considered, we have
\begin{flalign}
&\iint_{D_0} \left(1+{\textstyle \frac{\nu\sin^2\alpha(s)}{1-\nu}}\right)\frac{{\rm sgn}(P)}{2\pi d_P} K(\xi,\eta,t)d\xi d\eta \nonumber\\
=&\frac{1}{2\pi }\left(1+{\textstyle \frac{\nu\sin^2\alpha}{1-\nu}}\right)
\frac{d}{\left(\frac{1}{2}\left(1+{\textstyle \frac{\nu\sin^2\alpha}{1-\nu}}\right) t\right)^2+d^2}+C_1^{\kappa,\alpha,d/t,t}+O(t^{\frac{1}{5}})
+O\left(\frac{t}{r_1^2}\right).\label{eqn:1/d-g}
\end{flalign}

Summarizing the contributions of Case 1 and the three orders in Eqs.~\eqref{eqn:0-g},  \eqref{eqn:log-d-g}, and \eqref{eqn:1/d-g} of Case 2, for small $t$, we have
\begin{flalign}
u_t(0,d,t)=&-\frac{1}{2\pi }\left(1+{\textstyle \frac{\nu\sin^2\alpha}{1-\nu}}\right)\frac{d}{\left(\frac{1}{2}\left(1+{\textstyle \frac{\nu\sin^2\alpha}{1-\nu}}\right) t\right)^2+d^2}-C_1^{\kappa,\alpha,d/t,t} \nonumber\\
& +\frac{1+\nu(1-3\sin^2\alpha)}{4\pi(1-\nu)}\kappa \log\left[
{\textstyle \left(1+\frac{\nu}{1-\nu}\sin^2\alpha\right)} C_2\sqrt{4d^2+t^2}\right]-A(0,0)\nonumber\\
&+O\left(\frac{t}{r_1^2}\log r_1\right) +O(r_1\log t)
+O\left(t^{\frac{1}{3}}\log t\right).
\end{flalign}
Recall that here $A(0,0)$ is the $O(1)$ contribution in Eq.~\eqref{eqn:I00-asymp} at the point $P=(0,d)$.
Integrate it over time $[0,t]$, we have
\begin{flalign}
u(0,d,t)=&\frac{1}{2}+\frac{1}{\pi}\mathrm{arctan}\frac{2d}{\left(1+{\textstyle \frac{\nu\sin^2\alpha}{1-\nu}}\right)t}
-\int_0^t C_1^{\kappa,\alpha,d/\tau,\tau}d\tau
\nonumber\\
&+\frac{1+\nu(1-3\sin^2\alpha)}{4\pi(1-\nu)}\kappa t\log \left[
{\textstyle \left(1+\frac{\nu}{1-\nu}\sin^2\alpha\right)}e^{-1} C_2\sqrt{4d^2+t^2}\right]+2d\arctan\frac{t}{2d}
 -A(0,0)t\nonumber\\
& +O\left(\frac{t^2}{r_1^2}\log r_1\right)+O(r_1t\log t)
+O\left(t^{\frac{4}{3}}\log t\right)\label{eqn:velocity0000}\\
=&\frac{1}{2}+\frac{2d}{\pi\left(1+{\textstyle \frac{\nu\sin^2\alpha}{1-\nu}}\right)t}
+\frac{1+\nu(1-3\sin^2\alpha)}{4\pi(1-\nu)}\kappa t\log t
-C_0t -A(0,0)t\nonumber\\
& +O\left(\frac{t}{r_1^2}\log r_1\right)+O(r_1\log t)
+O\left(t^{\frac{1}{3}}\log t\right)+O\left(\frac{d}{t}\right)+O\left(\frac{d^2}{t^2}\log t\right),\label{eqn:u-convolution}
\end{flalign}
where
\begin{flalign}
C_0=\frac{1}{t}\int_0^t C_1^{\kappa,\alpha,d/\tau,\tau}d\tau-\frac{1+\nu(1-3\sin^2\alpha)}{4\pi(1-\nu)}\kappa \log \left(
{\textstyle \left(1+\frac{\nu}{1-\nu}\sin^2\alpha\right)}e^{-1} C_2\right).
\end{flalign}

Letting $u(0,d,t)=\frac{1}{2} $, i.e., the dislocation travel a small distance $d$ within time $t$, the velocity of the dislocation is $v=\frac{d}{t}$, which gives
\begin{flalign}
v=& \frac{\pi}{2}\left(1+{\textstyle \frac{\nu\sin^2\alpha}{1-\nu}}\right)t\left[-\frac{1+\nu(1-3\sin^2\alpha)}{4\pi(1-\nu)}\kappa \log t+C_0+A(0,0)\right.\nonumber\\
 &\left.+O\left(\frac{t}{r_1^2}\log r_1\right)+O(r_1\log t)
+O\left(t^{\frac{1}{3}}\log t\right)+O\left(\frac{d}{t}\right)+O\left(\frac{d^2}{t^2}\log t\right)\right]. \label{eqn:velt}
\end{flalign}
Here $C_0$ is the $O(1)$ contribution to the velocity generated by dislocation core, which also appears in the available discrete dislocation dynamics methods \cite{GB1976,kubin1992,Ghoniem2000,xiang2003level,quek2006level,cai2006non,Arsenlis2007,Zhao2012}. The approximation of the dislocation velocity in {Eq.~\eqref{eqn:velt}} holds when all the error terms in it $\ll 1$ as $t\rightarrow 0$.
{ Especially, if we set $r_1=t^\frac{1}{3}$, the error terms inside the brackets are $O\left(t^{\frac{1}{3}}\log t\right)+O\left(\frac{d}{t}\right)+O\left(\frac{d^2}{t^2}\log t\right)$, for $d\ll \frac{t}{\sqrt{|\log t|}}$ as $t\rightarrow 0$.}

In the threshold dislocation dynamics method, the solution $u$ is evolved accurately for a time period $\Delta t$ and then  is adjusted by thresholding. Therefore, the velocity of the dislocation { when it travels a small distance $d\ll \frac{\Delta t}{\sqrt{|\log \Delta t|}}$} is
\begin{flalign}
v=& \frac{\pi}{2}\left(1+{\textstyle \frac{\nu\sin^2\alpha}{1-\nu}}\right)\Delta t\left[-\frac{1+\nu(1-3\sin^2\alpha)}{4\pi(1-\nu)}\kappa \log \Delta t+C_0+A(0,0)\right.\nonumber\\
 &{ \left.
+O\left(\Delta t^{\frac{1}{3}}\log \Delta t\right)+O\left(\frac{d}{\Delta t}\right)+O\left(\frac{d^2}{\Delta t^2}\log \Delta t\right)\right]. } \label{eqn:vel0000}
\end{flalign}

{ Recall that as explained in Sec.~\ref{sec:Thresh},  $\Delta t$  is the dislocation core parameter in the dimensionless form presented at the end of Sec.~\ref{sec:PN}. In a threshold dynamics method, it is necessary that the moving front goes across at least one grid point. This implies that Eq.~\eqref{eqn:vel0000} should hold for $d=O(\Delta x)$. In this case, the error in the dislocation velocity in  Eq.~\eqref{eqn:vel0000} is $O\left(\Delta t^{\frac{1}{3}}\log \Delta t\right)+O\left(\frac{\Delta x}{\Delta t}\right)+O\left(\frac{\Delta x^2}{\Delta t^2}\log \Delta t\right)$, which requires $\Delta x \ll \frac{\Delta t}{\sqrt{|\log \Delta t|}}$.
}


\subsection{Dislocation velocity with applied stress}

For the full threshold dislocation dynamics formulation in Eqs.~\eqref{eqn:1-1} and \eqref{eqn:1-2}, which taking into consider the motion of dislocations under the applied stress $\sigma^{\rm app}$, it can be calculated using the method in the previous subsection that the corresponding dislocation velocity is
\begin{flalign}
v=& \frac{\pi \Delta t}{2}\left(1+\frac{\nu\sin^2\alpha}{1-\nu}\right)\left(-\frac{1+\nu(1-3\sin^2\alpha)}{4\pi(1-\nu)}\kappa \log \Delta t+C_0+A(0,0)+\sigma^{\rm app}\right), \label{eqn:vel-full}
\end{flalign}
whose errors are the same as those in Eq.~\eqref{eqn:vel0000}. If $\pmb \tau\cdot\mathbf b<0$, there is a negative sign in this velocity formulation.

A special case is that a straight dislocation moving under the applied stress $\sigma^{\rm app}$. In this case, the dislocation is moving only under the applied stress, and its velocity generated by the threshold dynamics method is
\begin{flalign}
v=& \frac{\pi \Delta t}{2}\left(1+ \frac{\nu\sin^2\alpha}{1-\nu}\right)\sigma^{\rm app}. \label{eqn:vel-straight}
\end{flalign}
Recall that $\alpha$ is the angle between the dislocation line direction $\pmb \tau$ and the Burgers vector $\mathbf b$. As above, if $\pmb \tau\cdot\mathbf b<0$, there is a negative sign in this velocity formulation.

The resulting dislocation velocity given in Eq.~\eqref{eqn:vel-full} agrees with those in the available discrete dislocation dynamics methods \cite{GB1976,kubin1992,Ghoniem2000,xiang2003level,quek2006level,cai2006non,Arsenlis2007,Zhao2012}, except for the
 orientation-dependent  prefactor $\frac{\pi \Delta t}{2}\left(1+\frac{\nu\sin^2\alpha}{1-\nu}\right)$. A numerical method to correct this specific anisotropic dislocation mobility caused by the dislocations stress kernel will be presented in the next section.

\section{Correcting dislocation velocity by spatial variable stretching}\label{sec:CTDD}

In our threshold dislocation dynamics method, both the evolution kernel and time step have physical meanings: the evolution kernel $K(x,y,t)$ in Eq.~\eqref{eqn:Kt} comes from the kernel for the elastic interaction between dislocations, and the time step $\Delta t$, in the dimensionless form of the equation, reflects the dislocation core radius.
These physical meanings  impose  restrictions in the implementation of the threshold dislocation dynamics method.

First, in order for the dislocation velocity to be accurate for two leading orders, in the dimensionless form of the equation, the time step $\Delta t$ has to be equal with the dislocation core radius, which has to be small due to the fact that the dislocation core size should be much less than the size of the domain. As a result, the effective dislocation velocity $v$, which is proportional to $\Delta t$ as shown by the formulation in Eq.~\eqref{eqn:vel-full}, is also small. The small dislocation velocity imposes a {severe} numerical limitation on the spatial grid constant, because the dislocation has to move across at least one spatial grid in the threshold dynamics method.

Moreover, due to the anisotropic dislocation stress kernel $K(x,y,t)$ in Eq.~\eqref{eqn:Kt}, the dislocation velocity obtained in Eq.~\eqref{eqn:vel-full} has an orientation-dependent mobility, i.e. the prefactor $\frac{\pi \Delta t}{2}\left(1+\frac{\nu\sin^2\alpha}{1-\nu}\right)$, where $\alpha$ is the angle between the dislocation line direction $\pmb \tau$ and its Burgers vector $\mathbf b$.
This anisotropic mobility is not necessarily the mobility of dislocations from physics.

In this section, we present a numerical method to solve the above two problems in the threshold dislocation dynamics. The numerical method is based on stretching of the spatial variables in the threshold dislocation dynamics equation in the convolution step. This method is able to correct dislocation mobility to any form. This stretching method is also able to enlarge the dislocation velocity, i.e.,
to speed up the dislocation motion, so that the threshold dislocation dynamics methods still applies when the numerical grid constant is not that small.

%

  \begin{figure}[!htb]
	\centering
	\includegraphics[width=4.0in]{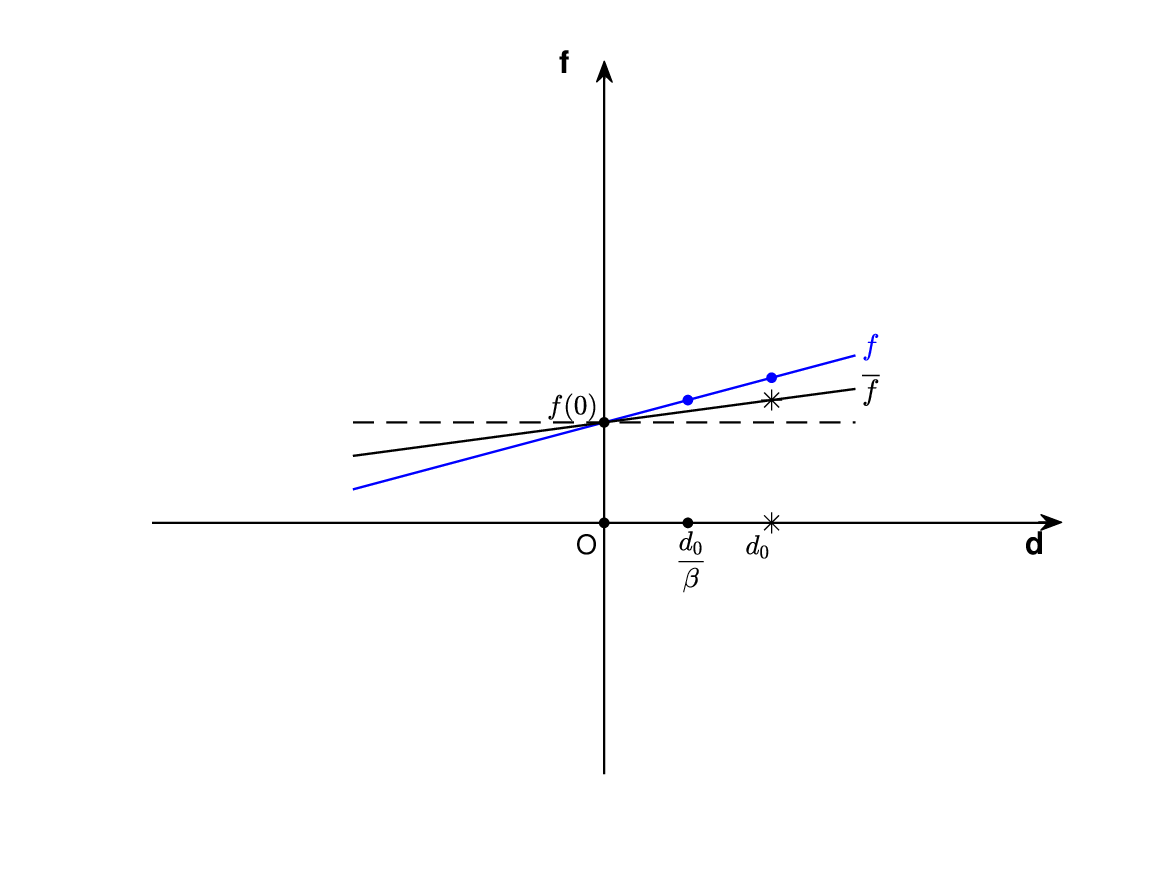}
	\caption{For a linear function $f(d)$,
  stretching in variable $d$ by a factor $\beta$ is equivalent to  stretching in $f$ (after subtraction of $f(0)$) by the factor $1/\beta$. The stretched function is $\overline{f}$. }
	\label{fig:1dctdd}
\end{figure}

  Here we illustrate the main idea of the variable stretching method. For a one dimensional linear function $f(d)$,
  stretching in variable $d$ by a factor $\beta$ is equivalent to  stretching in $f$ by the factor $1/\beta$, as shown in Fig.~\ref{fig:1dctdd}.
  Denoting $\overline{f}$ to be the stretched linear function, at a point $d_0$, we have
\begin{equation}
\overline{f}(d_0)=f\left(\frac{d_0}{\beta}\right)=f(0)+\frac{1}{\beta}[f(d_0)-f(0)].\label{eqn:stretch-linear}
\end{equation}
In the threshold dislocation dynamics method, near the dislocation, the solution $u$ before thresholding is approximately a linear function of the signed distance $d$ to the dislocation.
Details of the numerical methods based on variable stretching to correct the mobility and to rescale the velocity in the threshold dislocation dynamics method will be presented in the following two subsections.

\subsection{Correcting the mobility}
In this subsection, we present the variable stretching method to correct the dislocation mobility.  The method will be demonstrated based on the case of isotropic mobility, i.e., eliminating the anisotropic factor $1+\frac{\nu\sin^2\alpha}{1-\nu} $ that appears in the velocity formula in Eq.~\eqref{eqn:vel-full}.

\begin{figure}[htbp]
	\centering
	\includegraphics[width=3.0in]{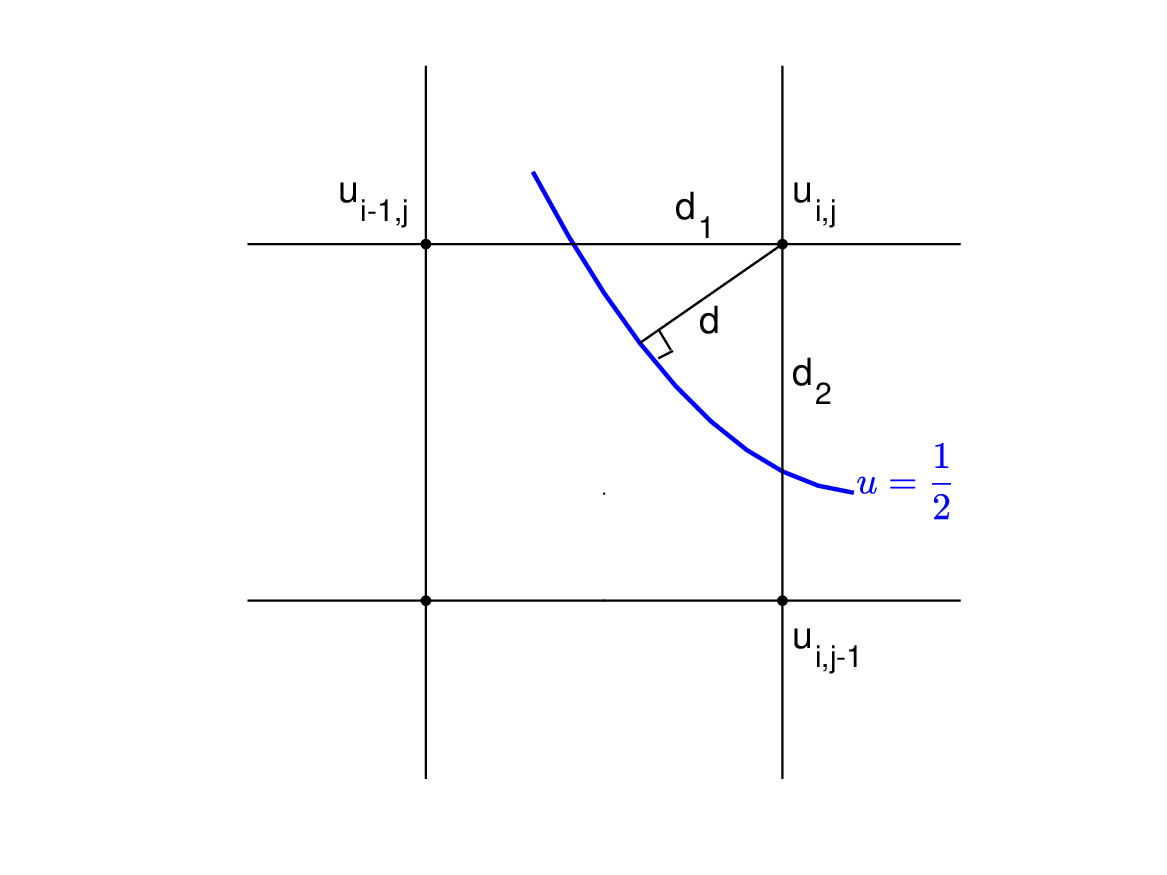}
	\caption{Solution $u_{ij}$ on the discrete grid points and the distance from a grid point to the dislocation. }
	\label{fig:2dctdd}
\end{figure}

Consider the solution $u(x,y,t)$ on the discrete grid points after a convolution step, as
 shown in Fig.~\ref{fig:2dctdd}, where the dislocation is given by $u=\frac{1}{2}$.
 Let $u_{ij}$ be the value of $u$ at discrete grid point $(x_i,y_j)$, and $d$ is the distance from the point $(x_i, y_j)$ to the dislocation.
 It can be calculated that
\begin{equation}
	d=\frac{d_1d_2}{\sqrt{d_1^2+d_2^2}},\label{eqn:d}
\end{equation}
where
\begin{equation}
	d_1=\frac{u_{i,j}-0.5}{u_{i,j}-u_{i-1,j}}\Delta x, \ \ \
	d_2= \frac{u_{i,j}-0.5}{u_{i,j}-u_{i,j-1}}\Delta y.
\end{equation}
The angle $\alpha$ between the dislocation line direction and the Burgers vector $\mathbf b$, which is in the $x$ direction here, can be calculated by
\begin{equation}
\alpha=\arctan\frac{d_2}{d_1}. \label{eqn:alpha}
\end{equation}

The idea of correcting dislocation mobility
  is to stretch u by stretching the space coordinate perpendicular to the dislocation direction with a factor $1/ (1+\frac{\nu\sin^2\alpha}{1-\nu})$ to eliminate the anisotropic coefficient $1+\frac{\nu\sin^2\alpha}{1-\nu} $.
  More precisely,
  from Eqs.~\eqref{eqn:u-convolution} and \eqref{eqn:vel0000}, and $v=d/\Delta t$, if $d$ is stretched to
  \begin{equation}
  \overline{d}=\frac{d}{1+\frac{\nu\sin^2\alpha}{1-\nu}},\label{eqn:stretch}
  \end{equation}
   the dislocation mobility will be isotropic, and the dislocation velocity will be corrected to
  \begin{flalign}
v=& \frac{\pi \Delta t}{2}\left(-\frac{1+\nu(1-3\sin^2\alpha)}{4\pi(1-\nu)}\kappa \log \Delta t+C_0+A(0,0)+\sigma^{\rm app}\right). \label{eqn:vel-full-adjusted}
\end{flalign}

Near the dislocation, the solution $u$ before thresholding is approximately a linear function of $d$. The stretching will following the stretching of one-dimensional linear function shown in Eq.~\eqref{eqn:stretch-linear}.


Assume that at time step $t_{n+1}$, the solution $u$ obtained after the convolution step is $\widetilde {u^{n+1}}$.
For the desired stretching in $d$ in Eq.~\eqref{eqn:stretch} right after the convolution step (before the thresholding step), using Eq.~\eqref{eqn:stretch-linear}, we have the formula for the stretched solution $\overline {u^{n+1}}$
\begin{equation}
\overline {u^{n+1}}= u_{\rm dis}+\left(1+\frac{\nu\sin^2 \alpha^{n}}{1-\nu}\right)(\widetilde {u^{n+1}}-u_{\rm dis}).\label{eqn:u-adjusted}
\end{equation}
Here $u_{\rm dis}$ is
\begin{equation}
u_{\rm dis}=\widetilde {u^{n+1}}-\frac{2}{\pi \Delta t}\frac{d^{n}}{1+\frac{\nu\sin^2\alpha^n}{1-\nu}}, \label{eqn:u-dis}
\end{equation}
which is the value independent of $d$ based on Eq.~\eqref{eqn:u-convolution}.

In summary, in the time step from $t_n$ to $t_{n+1}$, we compute $\widetilde {u^{n+1}}$ following Eq.~\eqref{eqn:formulation}, and $u_{\rm dis}$ by Eq.~\eqref{eqn:u-dis}, and then stretched solution $\overline {u^{n+1}}$ by Eq.~\eqref{eqn:u-adjusted}. We then perform thresholding for  $\overline {u^{n+1}}$ according to Eq.~\eqref{eqn:thresholding1} (or \eqref{eqn:thresholding2}) for the case of multiple dislocations). Following this procedure, we will have the desired dislocation velocity with isotropic mobility in Eq.~\eqref{eqn:vel-full-adjusted}.
Note that any physically meaningful anisotropic dislocation mobility can be assigned by this approach in the threshold dislocation dynamics method.

\subsection{Rescaling the velocity to a larger value}

We can further rescale the velocity to a larger value using the variable stretching algorithm shown above.

Note that in a threshold dynamics method, the moving front needs to move across at least one spatial grid before the thresholding step in order for the front eventually moves after the thresholding.
It can be seen that the effective dislocation velocity $v$ in Eq.~\eqref{eqn:vel-full} is proportional to $\Delta t$. Thus over a time step $\Delta t$, the travel distance of the dislocation is of order $(\Delta t^2)$, which requires a very small spatial grid constant $\Delta x$ for the dislocation to move after the thresholding, and even smaller $\Delta x$ for the velocity to be accurate.
By stretching the spatial domain in the direction normal to the dislocation by a factor $\beta>1$,  the velocity is rescaled by a factor of $\beta$. In this way, we effectively accelerate the motion of the dislocation, allowing it to move with a $\Delta x$ that is not that small  and also reducing the error in the dislocation velocity.

Incorporating the above velocity rescaling, the numerical method by variable stretching to adjust the dislocation velocity is
\begin{equation}
\overline {u^{n+1}}= u_{\rm dis}+\frac{1}{\pi}\mathrm {arctan}\left(\frac{1}{\beta}\left(1+\frac{\nu\sin^2 \alpha^{n}}{1-\nu}\right)
\tan\left(\pi(\widetilde {u^{n+1}}-u_{\rm dis})\right)\right),\label{eqn:u-adjusted2}
\end{equation}
where the velocity rescaling factor $\beta\geq 1$. Here the $\tan(\cdot)$ and $\mathrm {arctan}(\cdot)$ functions are used to avoid large values of the stress, e.g., when two dislocations are very close to each other; see the velocity formula before linearization in Eq.~\eqref{eqn:velocity0000}.

In the numerical implementation, the effect of applied stress can be added after the velocity corrections. That is,
\begin{equation}
\overline {u^{n+1}}= u_{\rm dis}+\frac{1}{\pi}\mathrm {arctan}\left(\frac{1}{\beta}\left(1+\frac{\nu\sin^2 \alpha^{n}}{1-\nu}\right)
\tan\left(\pi(\widetilde {u^{n+1}}-u_{\rm dis})\right)\right)-\sigma^{\rm app}\Delta t,\label{eqn:u-adjusted3}
\end{equation}
where the initial solution at time $t_{n+1}$, $\widetilde {u^{n+1}}=K_{\Delta t}*1_{S_n}$, only comes from the convolution and does not include the effect of the applied stress.

Numerically, velocity rescaled by a factor of $\beta$ implies that the effective time step is increased to $\beta \Delta t$.

\section{Algorithm of threshold dislocation dynamics method}\label{sec:algorithm}

In this section, we summarize the algorithm of the threshold dislocation dynamics method. In addition to the standard two steps of convolution and thresholding, it also includes a step of  correction of dislocation mobility and rescaling of dislocation velocity.


\vspace{0.1in}
\noindent
\underline{\bf Algorithm of Threshold Dislocation Dynamics Method (TDMM)}
\vspace{0.05in}

1. Give the initial condition $u^0$. Set the time step $\Delta t$ which corresponds to the dislocation core radius. Choose the velocity rescaling factor $\beta\geq 1$. The effective time step is $\beta \Delta t$.

2. Evolve the solution $u$ from $t_n$ to $t_{n+1}$ without the applied stress:
\begin{equation*}
\widetilde {u^{n+1}}= K_{\Delta t}*1_{S_n}.
\end{equation*}

3. Find the value of the solution $\widetilde {u^{n+1}}$ at the location of the dislocation at time $t^n$:
\begin{equation*}
u_{\rm dis}=\widetilde {u^{n+1}}-\frac{2}{\pi \Delta t}\widetilde{d^{n}}=\widetilde {u^{n+1}}-\frac{2}{\pi \Delta t}\frac{d^{n}}{1+\frac{\nu\sin^2\alpha^n}{1-\nu}}.
\end{equation*}

4. Adjust the solution at time $t_{n+1}$ to correct the mobility and to rescale the velocity by factor $\beta$:
\begin{equation*}
\overline {u^{n+1}}= u_{\rm dis}+\frac{1}{\pi}\mathrm {arctan}\left(\frac{1}{\beta}\left(1+\frac{\nu\sin^2 \alpha^{n}}{1-\nu}\right)
\tan\left(\pi(\widetilde {u^{n+1}}-u_{\rm dis})\right)\right).
\end{equation*}

5. Add the effect of the applied stress:
\begin{equation*}
\overline {u^{n+1}}= \overline {u^{n+1}}-\sigma^{\rm app}\Delta t.
\end{equation*}

6. Find the distance $d^{n+1}$ to the dislocation (where $u= 1/2$) and the angle $\alpha^{n+1}$ between the dislocation line direction and Burgers vector at $t_{n+1}$, using $\overline {u^{n+1}}$ and following Eqs.~\eqref{eqn:d}-\eqref{eqn:alpha}.

7. Update the solution $u$ at $t_{n+1}$ using threshold:
\begin{equation*}
u^{n+1}= j, \ \ {\rm if}\ j-0.5<\overline {u^{n+1}}\leq j+0.5.
\end{equation*}

8. Repeat steps 2-7.

%
%
%
%
%
%
%

\section{Numerical Simulations} \label{sec:Num}

In this section, we perform some numerical simulations using the developed threshold dislocation dynamics method. The simulation domain is $[-\pi,\pi]\times[-\pi,\pi] $.
The simulation domain corresponds to a physical domain of size $157b\times 157b$. That is, $b=2\pi/157\approx0.04$. The Burgers vector of the dislocations is $\mathbf b=(b,0)$.  We choose time step $\Delta t=0.16$, { meaning that the dislocation core radius is $0.16=4b$.} As specified in Sec.~\ref{sec:PN}, the length unit of the simulation domain is $l_0$, the time unit is $l_0/M_p\mu$, and the stress unit is $\mu b/l_0$. The Poisson ratio $\nu=1/3$ unless otherwise specified. {The factor $\beta=1$ unless otherwise specified.} Simulation results will be compared  with those of theoretic predictions \cite{anderson2017theory} and discrete dislocation dynamics simulations \cite{xiang2003level,xiang2004level}.

\subsection{Motion of a straight edge dislocation under applied stress}

\begin{figure}[htbp]
	\centering
	{\includegraphics[width=2.2in]{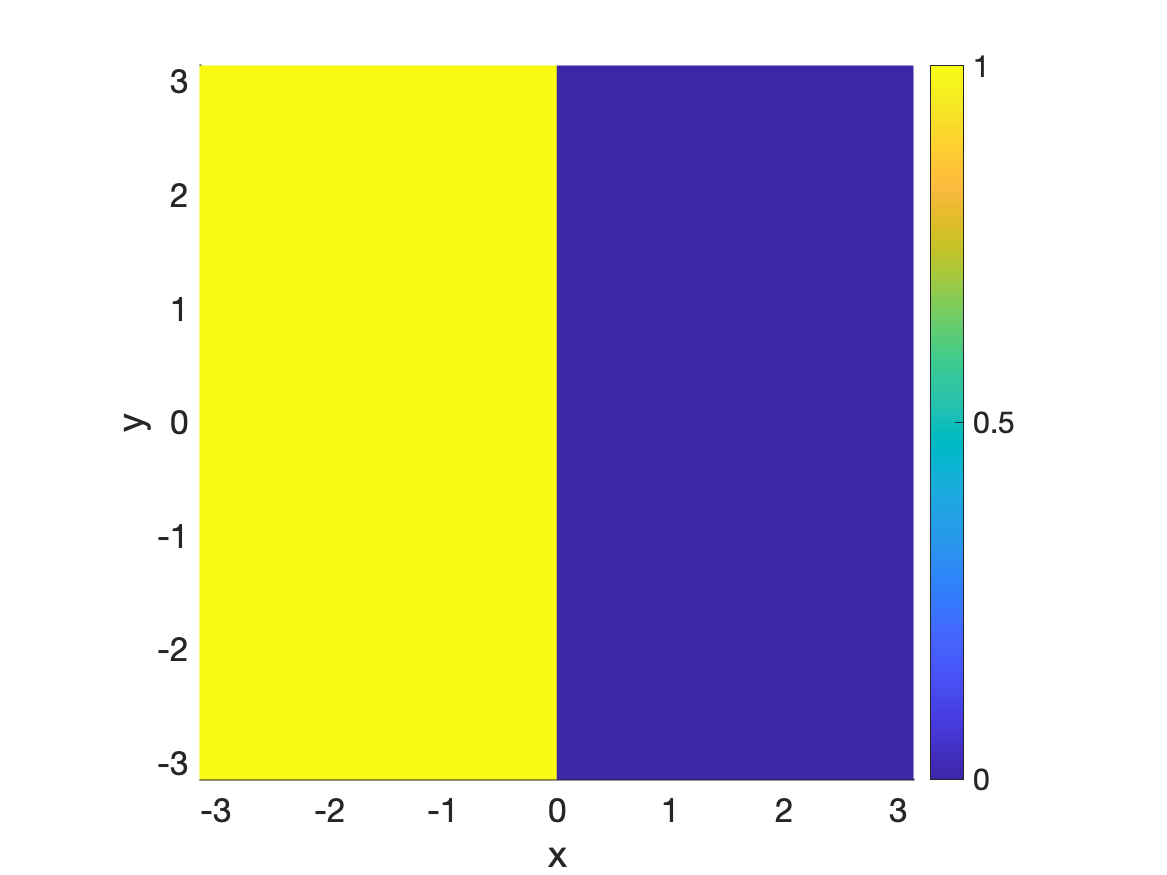}}
	\caption{Initially, an edge dislocation is located at $x=0$, with $u=1$ when $x>0$ and $u=0$ when $x<0$. The Burgers vector of the edge dislocation is $\mathbf b=(b,0)$.}
	\label{fig:straight_u}
\end{figure}

In this subsection, we simulate the motion of a straight edge dislocation under applied stress.
The dislocation is parallel to the $y$ axis and in the $+y$ direction.
Initially, the dislocation is located at $x=0$, and the initial condition of $u$ is shown in Fig.~\ref{fig:straight_u}.
The dislocation will move to the left with an applied stress $\sigma^{\rm app}=\sigma>0$.

As obtained in Eq.~\eqref{eqn:vel-full}, using the threshold dislocation dynamics method without velocity correction, the theoretic value of the dislocation velocity in this case is
\begin{equation}
		v_{a}=\frac{\pi\Delta t }{2(1-\nu)}\sigma.\label{eqn:edge-v0}
\end{equation}
For this edge dislocation, the anisotropic coefficient is $1+\frac{\nu\sin^2\alpha}{1-\nu}=1+\frac{\nu}{1-\nu}=\frac{1}{1-\nu} $. After velocity correction, this anisotropic coefficient is converted to $1$, and the theoretic value of the dislocation velocity is
\begin{equation}
		v=\frac{\pi\Delta t }{2}\sigma.\label{eqn:edge-v}
\end{equation}

\begin{figure}[!htb]
	\centering
{\includegraphics[width=3.4in]{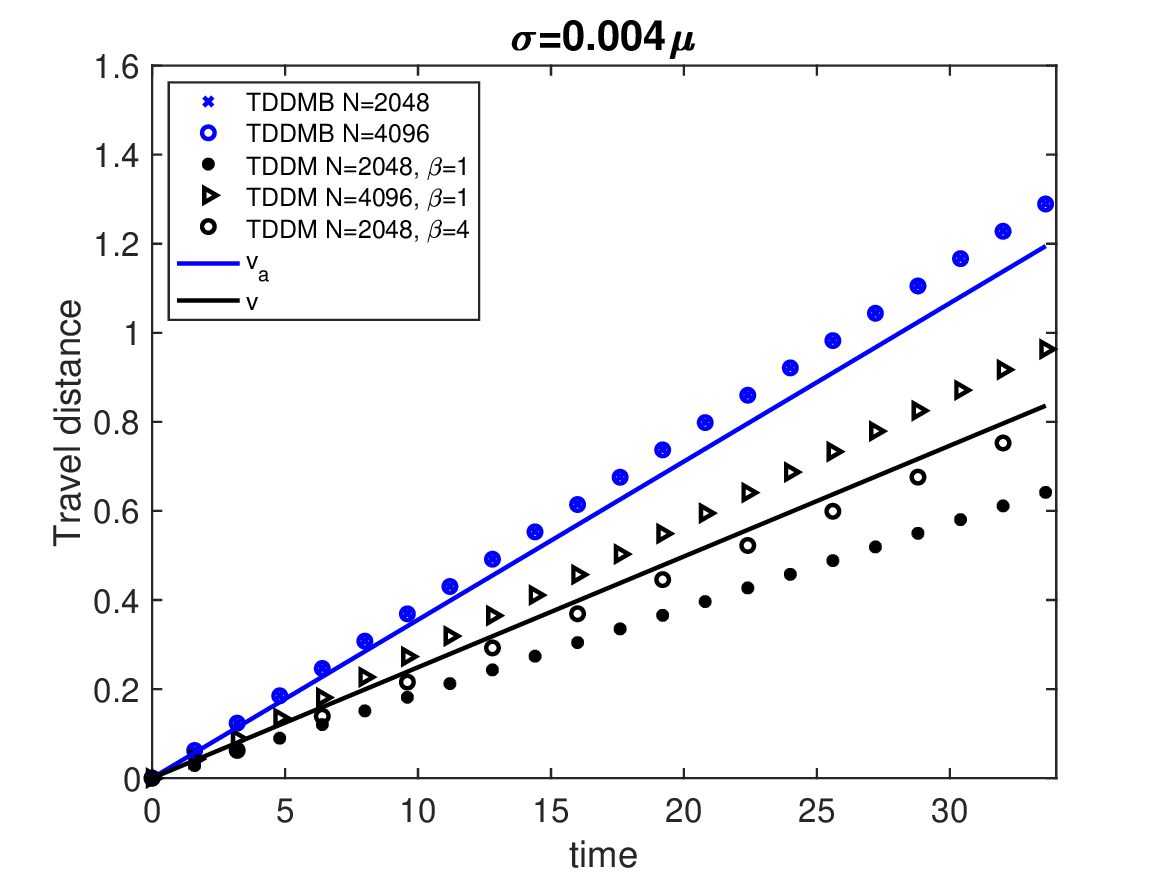}}

	\caption{Simulation results of motion of an edge dislocation under applied stress $\sigma=0.004\mu$, using the threshold dislocation dynamics method with mobility correction (TDDM) and without mobility correction (TDDMB).
		We compared the velocity using different rescaling factors, specifically $\beta=1$ and $\beta=4$, and different grid sizes, namely $N=2048$ and $N=4096$.
		The simulation results are compared with the theoretic values of the velocity with anisotropic mobility ($v_a$ in Eq.~\eqref{eqn:edge-v0}) and with the desired isotropic mobility ($v$ in Eq.~\eqref{eqn:edge-v}). The Poisson ratio $\nu=1/3$.
	}\label{fig:onethresh}
\end{figure}

Simulation results of the motion of this edge dislocation under applied stress  $\sigma=0.004\mu$ ($\sigma=0.1$ in the dimensionless form) under different numerical settings are shown in Fig.~\ref{fig:onethresh}.
As shown in { Fig.~\ref{fig:onethresh}}, when the numerical grid of the spatial domain is $2048\times 2048$ ($N=2048$),
 without correction of the dislocation mobility, i.e. $\beta=1$, the dislocation is approximately moving with the theoretical velocity $v_a$ in Eq.~\eqref{eqn:edge-v0} with anisotropic factor, whereas it has a systematic non-negligible difference compared with the desired isotropic velocity $v$ in Eq.~\eqref{eqn:edge-v};
  after numerical correction of the mobility, the error between the simulated dislocation velocity (the black dots) and the desired velocity (black line) is significantly reduced compared with the difference between the uncorrected simulated velocity (blue dots) and the desired velocity (black lines). When the velocity rescaling factor $\beta=4>1$, the error between the simulated dislocation velocity ({the black circles}) and the desired velocity (black line) is much smaller.

 We also perform simulation of the motion of this edge dislocation with a finer mesh $4096\times 4096$ ($N=4096$) and without velocity correction ($\beta=1$), and the results are shown in { Fig.~\ref{fig:onethresh}}. The error between the simulated dislocation velocity ({the black triangles}) and the desired velocity (black line) is also reduced compared with the results of a coarser mesh $N=2048$ {(black dots)}. Compared with the results {using a greater rescaling factor $\beta$ (black circles)}, it can be seen that in order to reduce the error, using a velocity rescaling factor $\beta>1$ is more effective than mesh refinement.

%


%


Simulation results of the motion of this edge dislocation under different values of the applied stress, with numerical correction of mobility and different numerical rescaling factors of velocity  are shown in Tables \ref{tab:table-1} and  \ref{tab:table-2}, in which the numerical grids of the spatial domain are $2048\times 2048$ ($N=2048$) and $4096\times 4096$ ($N=4096$), respectively.
These simulation results
demonstrate quantitatively that the numerical velocity of the dislocation is accelerated by a factor of $\beta>1$ when the velocity rescaling algorithm is applied. First consider the results for $N=2048$ shown in Table \ref{tab:table-1}.
When the applied stress is $\sigma=0.004\mu$,  without velocity rescaling, i.e. $\beta=1$, the relative error is about $23\%$; see also {Fig.~\ref{fig:onethresh}}. This error is reduced to $3.75\%$ when the velocity rescaling factor $\beta=4$ (see also {Fig.~\ref{fig:onethresh}}) and $0.1\%$ when $\beta=10$. Under a larger applied stress, these errors become smaller; this is because the dislocation is able to travel over more spatial grids within a time step. Under a smaller applied stress $\sigma=0.0008\mu$, the dislocation is not able to move without velocity rescaling. In fact, in this case, the dislocation is not able to move across a spatial grid with the given time step $\Delta t$. Using the velocity rescaling factor $\beta=4$ and $10$, the dislocation is able to move under this small applied stress, with errors in velocity about $10\%$. Using a smaller spatial grid of $N=4096$,  as shown in Table \ref{tab:table-2}, the simulations give more accurate values of dislocation velocity compared with those in the case of $N=2048$.

%

\begin{table}[!htbp]
	\caption{\label{tab:table-1}
		Velocity $\bar{v}$ of an edge dislocation obtained in simulations with different values of applied stress and rescaling factor $\beta$. $\Delta t=0.16$. $N=2048$. The Poisson ratio $\nu=1/3$.}
	\centering
	\begin{tabular}{|c|c|c|c|c|}
		\hline
		Applied stress  & $v$ &$\bar{v}$ \ with \ $\beta=1$&$\bar{v}$ \ with \   $\beta=4$ &$\bar{v}$ \ with \   $\beta=10$ \\
		\hline
		$0.0008\mu$ &0.0050 & 0.0=0$v$& 0.0192=3.84$v$ & 0.0575=11.5$v$    \\
		\hline
		$0.002\mu$ &0.0124 & 0.0197=1.55$v$& 0.0575=4.63$v$ & 0.1150=9.27$v$    \\
		\hline
		$0.004\mu$ &0.0249 &0.0197=0.77$v$ &0.0959=3.85$v$  & 0.2493=10.01$v$     \\
	    \hline
		$0.006\mu$ & 0.0374&0.0383=1.02$v$ & 0.1534=4.10$v$ & 0.3643=9.74$v$  \\
		\hline
		$0.008\mu$ & 0.0499 &0.0575=1.15$v$ &0.1917=3.84$v$  & 0.4985 =9.99$v$  \\
		\hline
	\end{tabular}
\end{table}

\begin{table}[htbp]
	\caption{\label{tab:table-2}   Velocity $\bar{v}$ of an edge dislocation obtained in simulations with different values of applied stress and rescaling factor $\beta$. $\Delta t=0.16$.      $N=4096$. The Poisson ratio $\nu=1/3$.}
	\centering
	\begin{tabular}{|c|c|c|c|c|}
	\hline
	Applied stress  & $v$ & $\bar{v}$ \ with \  $\beta=1$& $\bar{v}$ \ with \  $\beta=4$ &$\bar{v}$ \ with \   $\beta=10$ \\
	\hline
	$0.0008\mu$ &0.0050 & 0.0096=1.92$v$& 0.0192=3.84$v$ & 0.0497=9.94$v$    \\
	\hline
	$0.002\mu$ &0.0124 & 0.0096=0.77$v$& 0.0479=3.86$v$ & 0.1246=10.05$v$    \\
	\hline
	$0.004\mu$ &0.0249 &0.0288=1.16$v$ &0.0959=3.85$v$  & 0.2493=10.01$v$     \\
	\hline
	$0.006\mu$ & 0.0374&0.0383=1.02$v$ & 0.1534=4.10$v$ & 0.3739=10.00$v$  \\
	\hline
	$0.008\mu$ & 0.0499 &0.0479=0.96$v$ &0.2013=4.03$v$  & 0.4985 =9.99$v$  \\
	\hline
\end{tabular}
\end{table}

In short, our simulation results show that
the accuracy of the threshold dislocation dynamics method  can be improved by using a velocity rescaling factor $\beta>1$ and a finer spatial grid, and the velocity rescaling method is able achieve more significant improvment. These results demonstrate the effectiveness of the velocity rescaling method for obtaining more accurate results using coarser grids.

\subsection{Shrinking of a  circular dislocation under self stress}
Consider the evolution of a circular dislocation loop with initial radius $r=52.3b$ ($r=\frac{2\pi}{3}$ in the dimensionless form), and its center is at $(0,0)$. The loop will shrink under its self-stress. The numerical grid of the spatial domain is $1024\times 1024$ ($N=1024$).

\begin{figure}[!htb]
	\centering
	\subfigure{\includegraphics[width=1.45in]{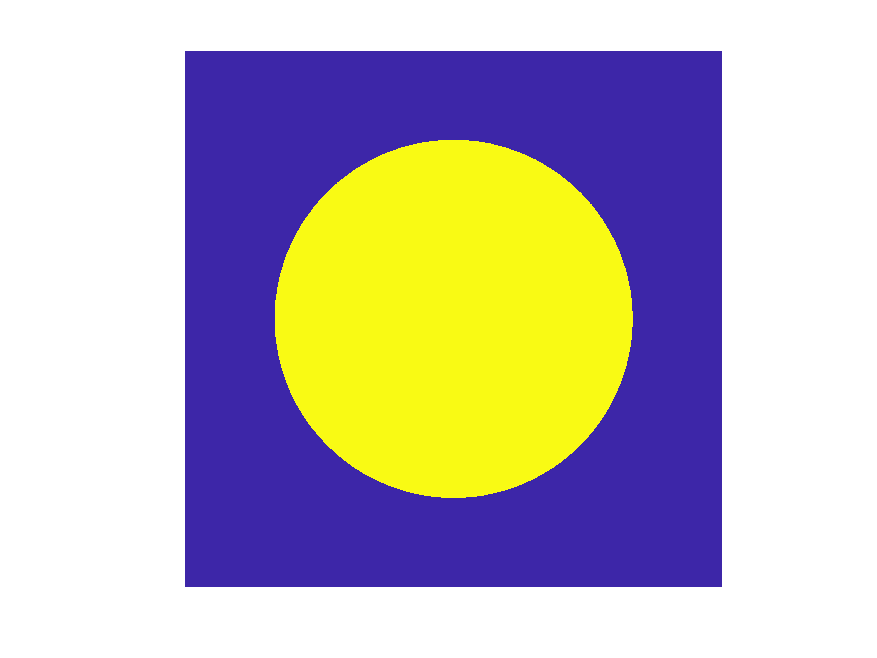}}\hspace{-0.2in}
	\subfigure{\includegraphics[width=1.45in]{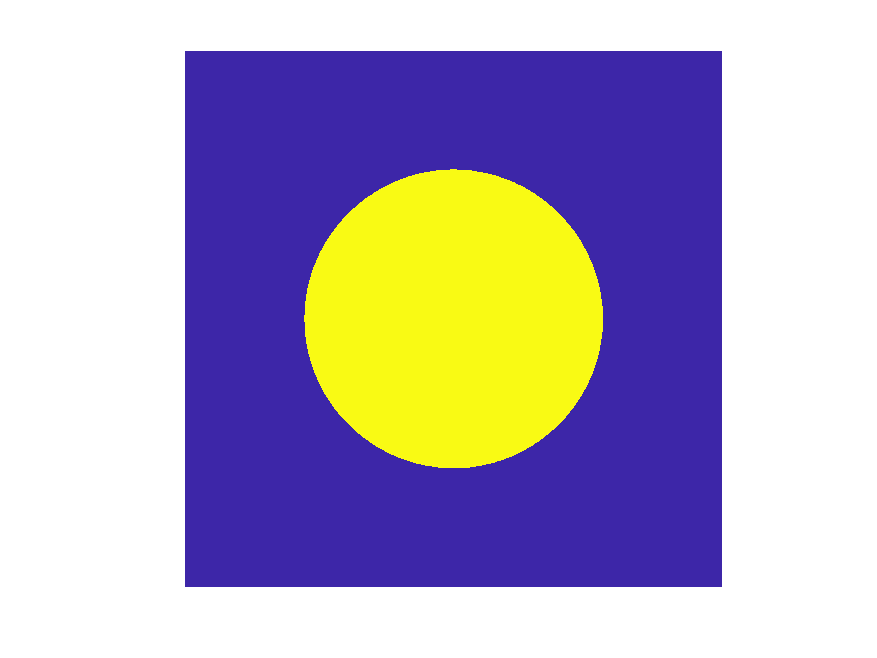}}\hspace{-0.2in}
	\subfigure{\includegraphics[width=1.45in]{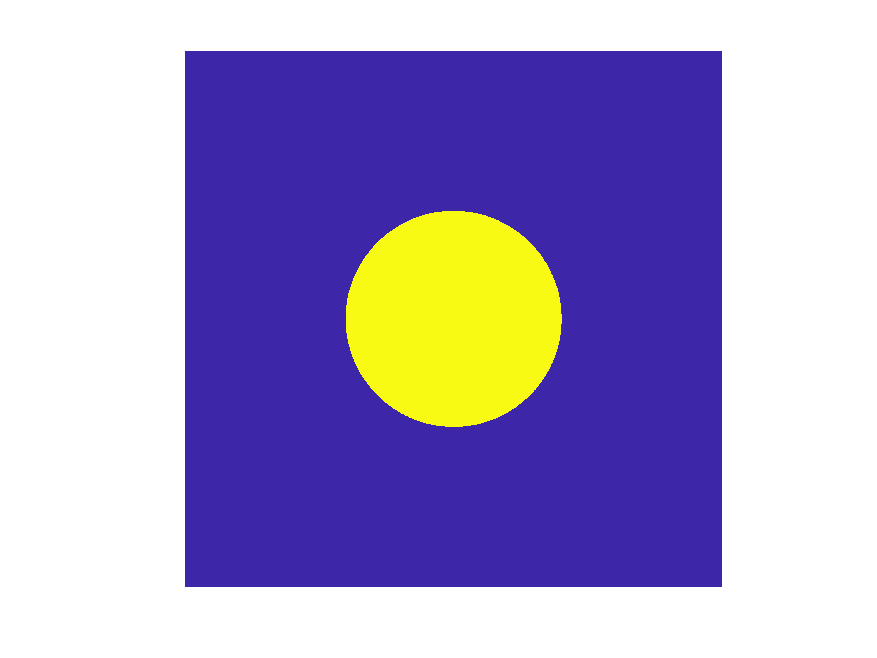}}\hspace{-0.2in}
	\subfigure{\includegraphics[width=1.45in]{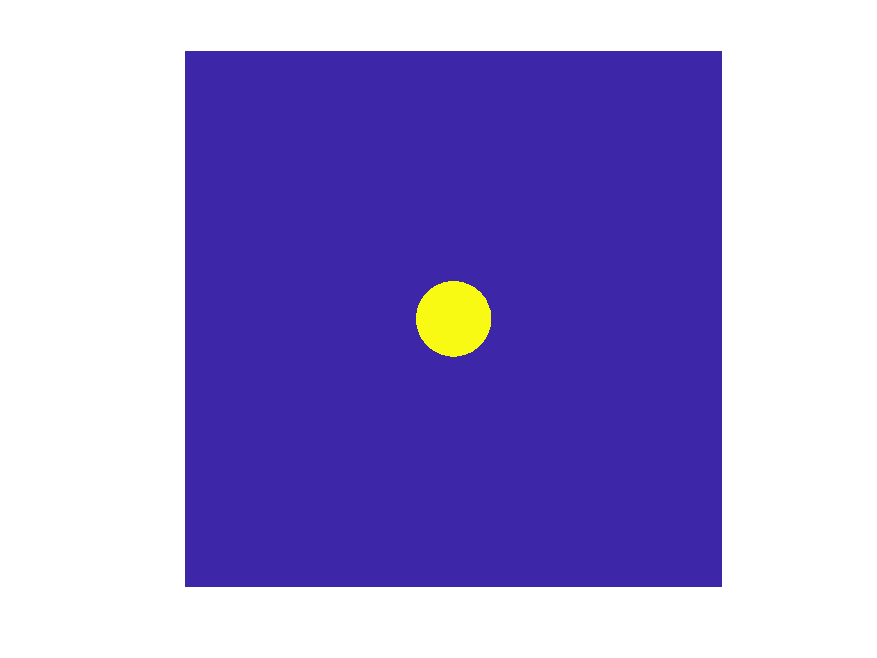}}
	\caption{Solution $u(x,y)$ for the evolution of an initially circular dislocation, where $\nu=0$. The yellow region is $u=1$ and the blue region is $u=0$, and the dislocation loop is the boundary between these two regions. }\label{fig:disu}
\end{figure}

\begin{figure}[htbp]
	\centering
	\subfigure{\includegraphics[width=2.7in]{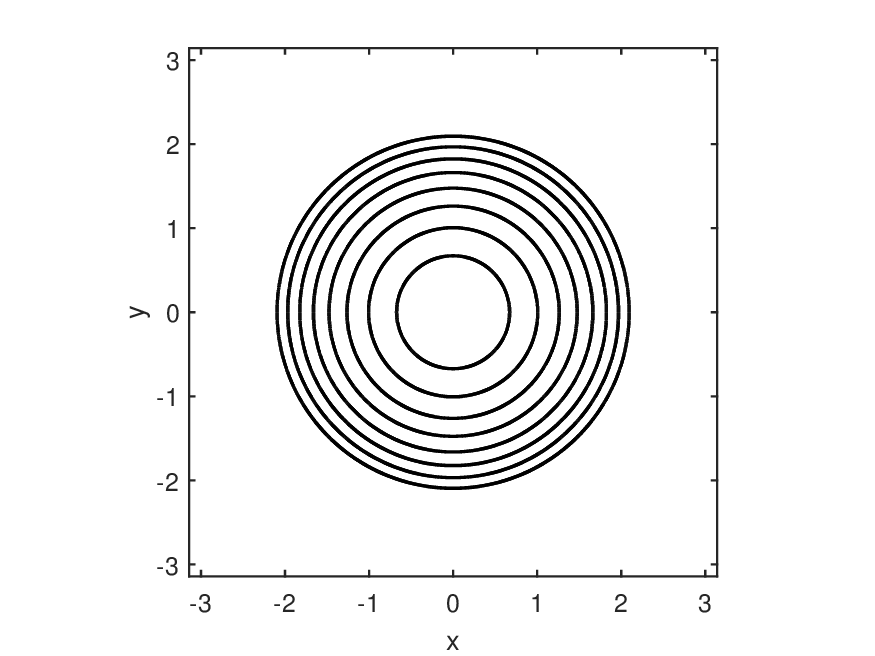}\label{fig:disnu0}}
  \subfigure{\includegraphics[width=2.7in]{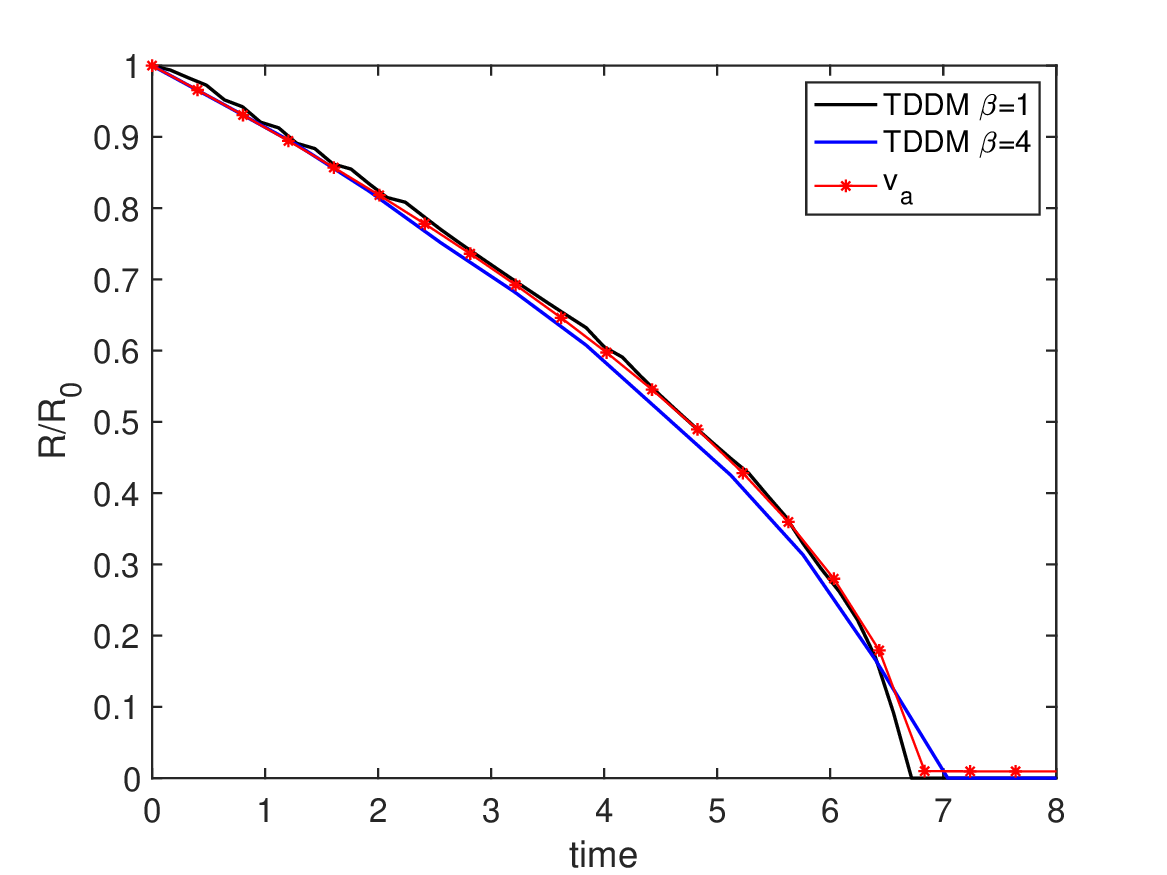}\label{fig:RChange}}
	\caption{(a) Simulation of an initial circular dislocation loop (the outermost circle)  shrinking under its self-stress when $\nu=0$. The dislocation loop is plotted at uniform time intervals starting with the outermost circle. The loop eventually disappears. (b) Radius $R$ of this circular loop during the evolution obtained using our threshold dislocation dynamics method (TDDM) with velocity rescaling factor $\beta=1$ and $\beta=4$, and comparison with the results given by the approximate velocity formula $v_a$ in \eqref{eqn:va-loop}. $R_0$ is the radius of the initial circular loop.}
\end{figure}

We firstly consider the case $\nu=0$. (This case is also corresponding to the physical process of  shrinking of a circular prismatic by climb, up to a factor  $1/(1-\nu)$ in velocity.) Fig.~\ref{fig:disu} shows the solution $u$ during the evolution. The evolution of dislocation (which is the boundary between the regions $u=0$ and $u=1$) is shown in Fig.~\ref{fig:disnu0}, and comparisons with the result using an approximate velocity formula are shown in Fig.~\ref{fig:RChange}.  The approximate velocity formula for the shrinking of this circular loop is
\begin{equation}
v_a=\frac{\Delta t}{8R}\log\frac{16R}{\Delta t}, \label{eqn:va-loop}
\end{equation}
which can be calculated from the general velocity formula in Eq.~\eqref{eqn:dis-v} before convolution with the kernel $\delta_\varepsilon$ (see, e.g. 4.2 in \cite{Gu}) and then averaging the velocity over a time interval of $\Delta t$. The asymptotic dislocation velocity in our threshold dislocation dynamics given in Eq.~\eqref{eqn:vel0000} agrees with this approximate velocity formula.
It can be seen from Fig.~\ref{fig:RChange} that the simulation results obtained by using our threshold dynamics method agree excellently with those by using the approximate velocity formula, which validates both methods. In this example, the velocity scaling factor $\beta=4>1$ plays a role of accelerating the simulation, i.e., effectively changing the time step from $\Delta t$ to $\beta \Delta t$.

\begin{figure}[htbp]
	\centering
	\subfigure{\includegraphics[width=3.0in]{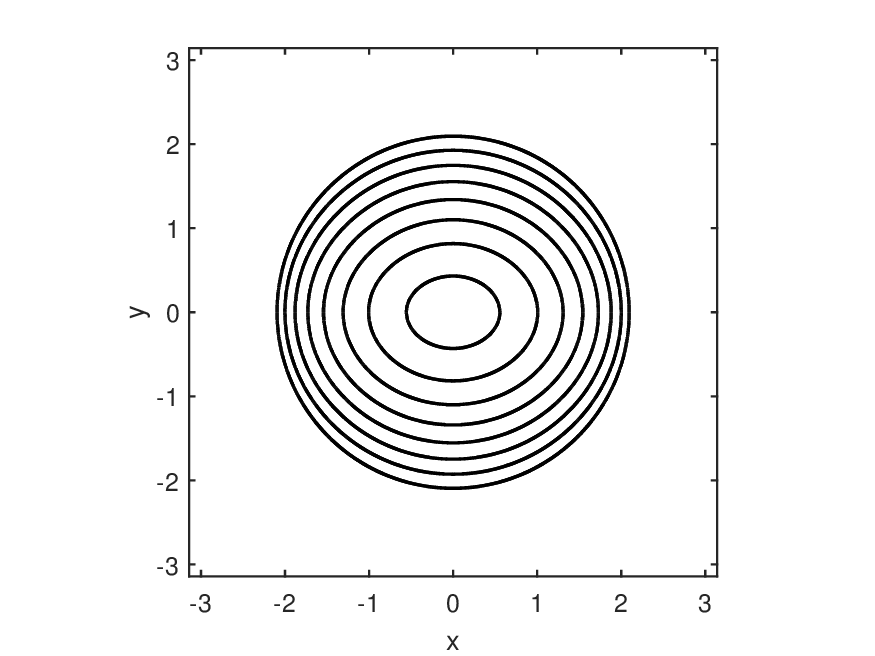}}
	\caption{Simulation of an initial circular dislocation loop (the outermost circle)  shrinking under its self-stress when $\nu=1/3$. The Burgers vector is in $+x$ direction. The initial circular dislocation loop shrinks and gradually becomes ellipse under its self-stress, and eventually disappears. The dislocation loop is plotted at uniform time intervals starting with the outermost circle.  }\label{fig:nvb}
\end{figure}

Simulation result for the evolution of this initially circular dislocation loop with Poisson ratio $\nu=1/3$ is shown
Fig.~\ref{fig:nvb}. Since the Burgers vector is in x direction, the dislocation is screw when the unit tangent vector is in $x$ direction and is edge when the unit tangent vector is in $y$ direction. The leading order shrinkage force is greater for screw dislocation segment than that on edge segments \cite{anderson2017theory,xiang2003level}. The asymptotic dislocation velocity in our threshold dislocation dynamics given in Eq.~\eqref{eqn:vel0000}, after corrected the anisotropic mobility factor $1+\frac{\nu\sin^2\alpha}{1-\nu}$ and neglecting the constant factor $\frac{\pi\Delta t}{2}$, is $\frac{1+\nu(1-3\sin^2\alpha)}{4\pi(1-\nu)}\kappa \log \frac{1}{\Delta t}$, which agrees with the known leading order velocity formula \cite{anderson2017theory,xiang2003level} (up to constant factor $\mu b^2$).
  The initially circular dislocation loop  gradually becomes ellipse when it is shrinking and eventually disappears.

\subsection{Evolution of two dislocation loops}

We simulate a system of two circular dislocation loops with same direction and radius $r=\frac{\pi}{3}$. The numerical grid of the spatial domain $N=1024$. The two dislocation loops will evolve under both their self stress and the stress due to the long-range interaction between them.

\begin{figure}[!htb]
	\centering
	\subfigure{\includegraphics[width=0.8in]{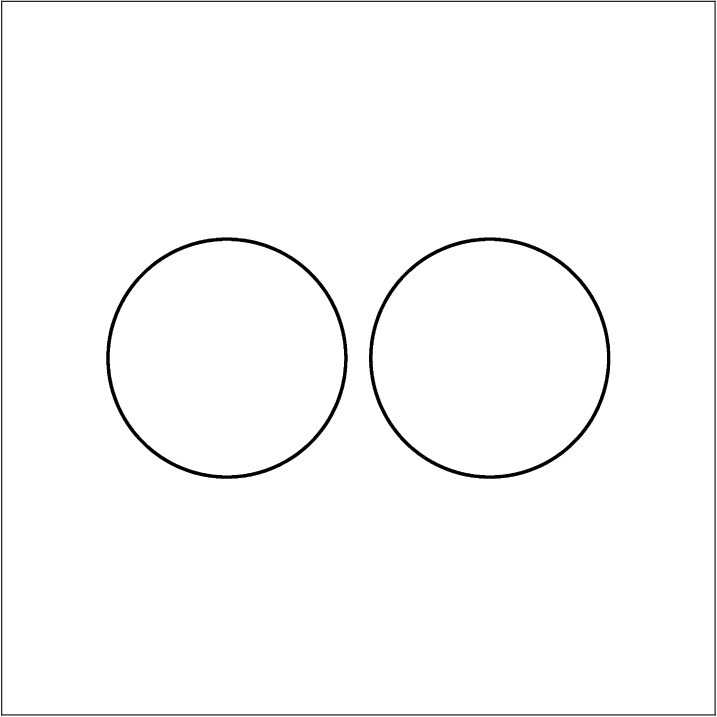}}\hspace{0.25in}
	\subfigure{\includegraphics[width=0.8in]{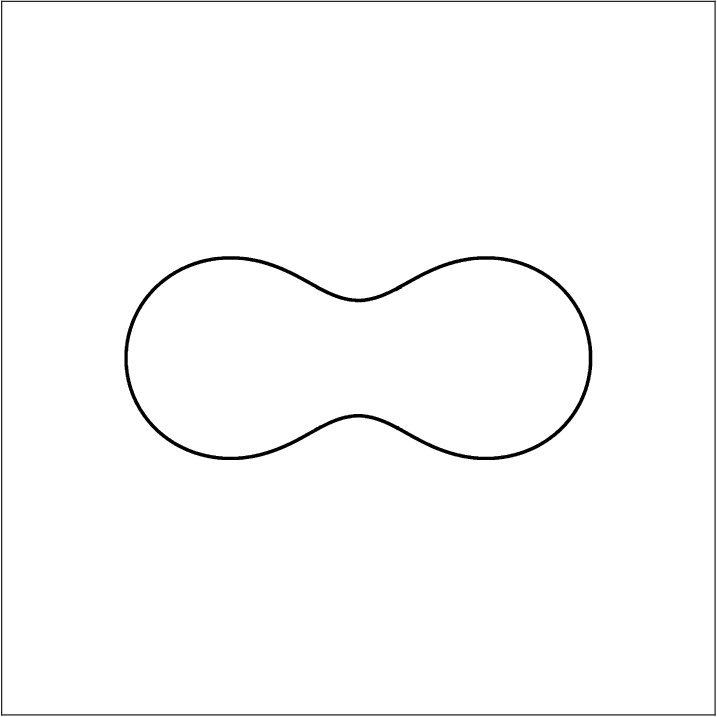}}\hspace{0.25in}
	\subfigure{\includegraphics[width=0.8in]{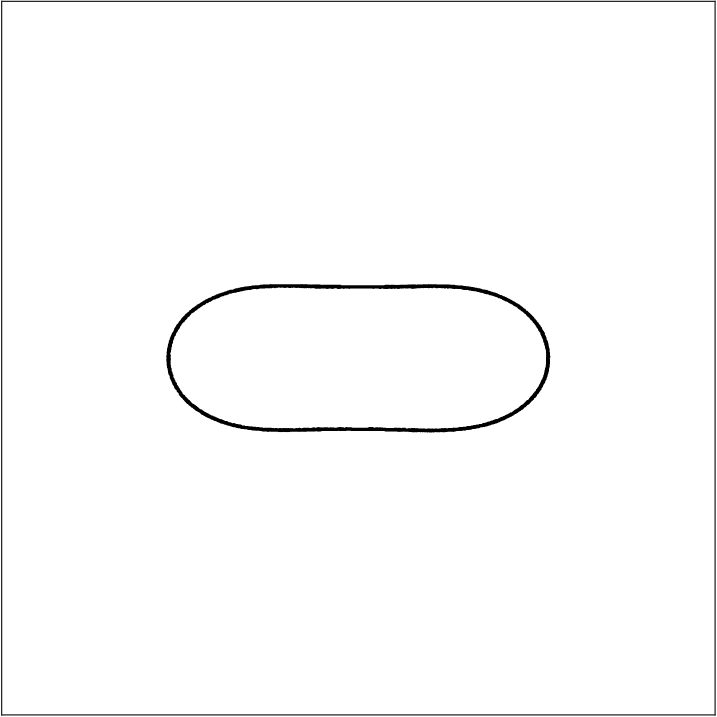}}\hspace{0.25in}
	\subfigure{\includegraphics[width=0.8in]{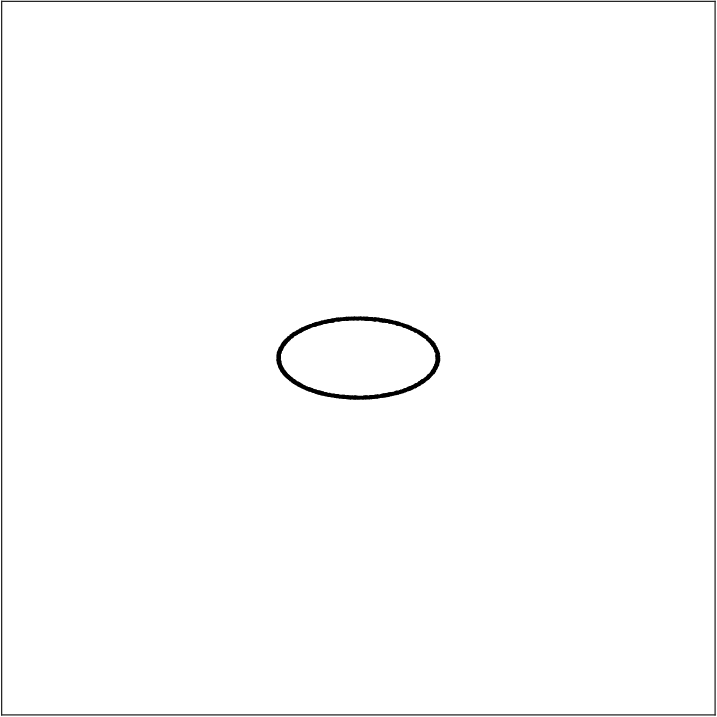}}\hspace{0.25in}
	\subfigure{\includegraphics[width=0.8in]{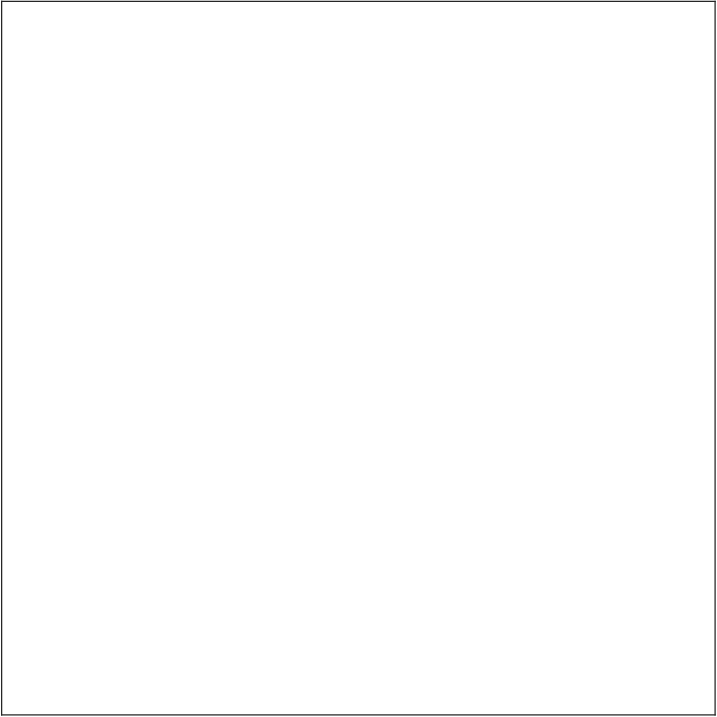}}
	\caption{Simulation results of evolution of two circular dislocation loops close to each other. The two loops are in the clockwise direction, and have the same Burgers vector that is in the horizontal direction. The Poisson ratio $\nu=1/3$.  }\label{fig:twocircon}
\end{figure}

We first simulate the evolution of the two circular dislocation loops in the clockwise direction and with their centers located at $(\frac{\pi}{3}+0.11,0)$ and $(-\frac{\pi}{3}-0.11,0)$, respectively; see Fig.~\ref{fig:twocircon}. These two loops initially are very close to each other. Due to the strong long-range interaction stress, the two loops first combine into a single loop, and then the single loop gets smoother as it shrinks under its self stress. The single loop eventually disappears.

\begin{figure}[!htb]
	\centering
	\subfigure{\includegraphics[width=0.8in]{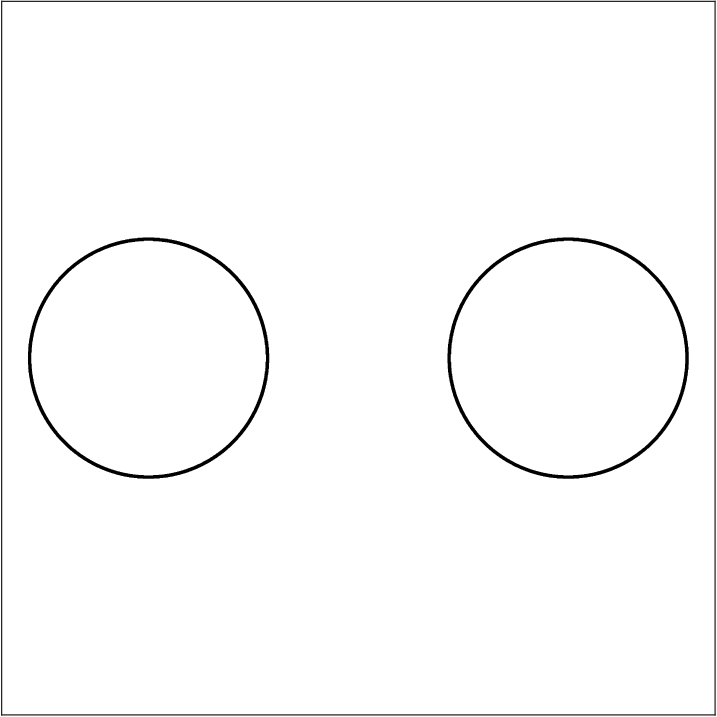}}\hspace{0.3in}
	\subfigure{\includegraphics[width=0.8in]{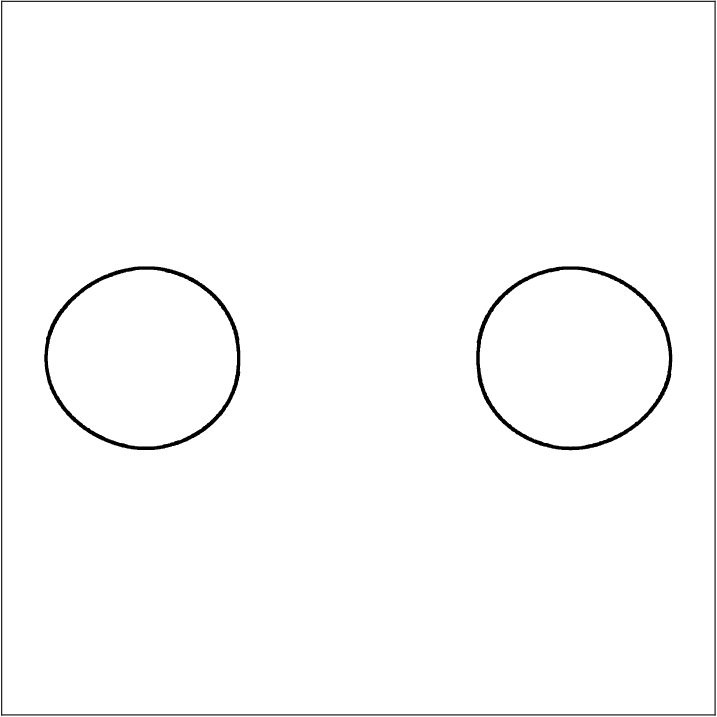}}\hspace{0.3in}
	\subfigure{\includegraphics[width=0.8in]{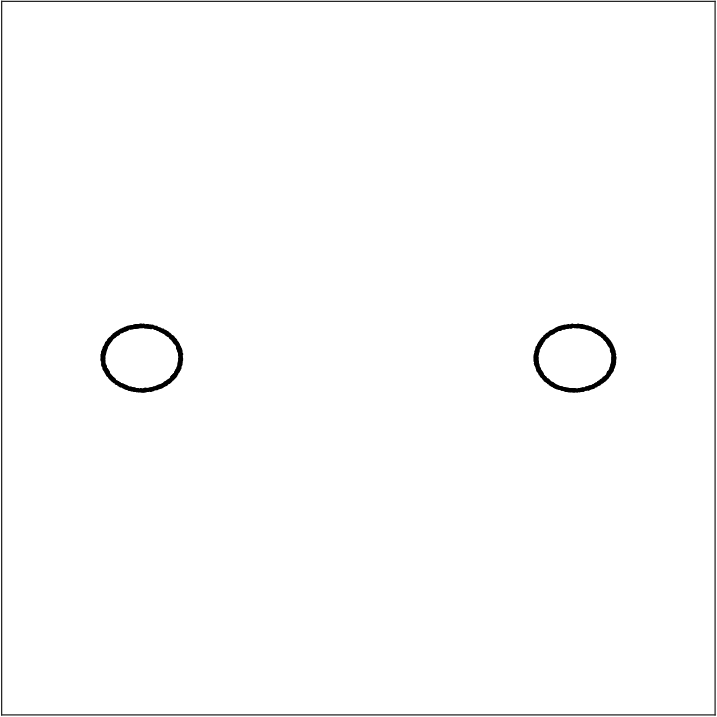}}\hspace{0.3in}
	\subfigure{\includegraphics[width=0.8in]{blank.eps}}
	\caption{Simulation results of evolution of two circular dislocation loops that are relatively separated initially. The two loops are in the clockwise direction, and have the same Burgers vector that is in the horizontal direction. The Poisson ratio $\nu=1/3$.     }\label{fig:twocirnon}
\end{figure}

We also simulate another case in which the initial two circular dislocation loops are relatively separated, with their centers located at $(\frac{\pi}{3}+0.8,0)$ and $(-\frac{\pi}{3}-0.8,0)$, respectively; see Fig.~\ref{fig:twocirnon} for the evolution of this system. Both loops are in the clockwise direction as in the previous case. In this case, the long-range interaction between the two loops is relatively weak compared with their self stress, and both loops are able to shrink under their self stress. The two loops eventually disappear.

These simulation results demonstrate that our threshold dislocation dynamics method can indeed correctly capture both the leading order ($O(\log\varepsilon)$, where $\varepsilon$ is the dislocation core width) curvature motion and the next order ($O(1)$) long-range interaction for the dynamics of dislocations.

\subsection{Dislocation bypassing particle }

\begin{figure}[!htb]
	\centering
	\subfigure{\includegraphics[width=1.0in]{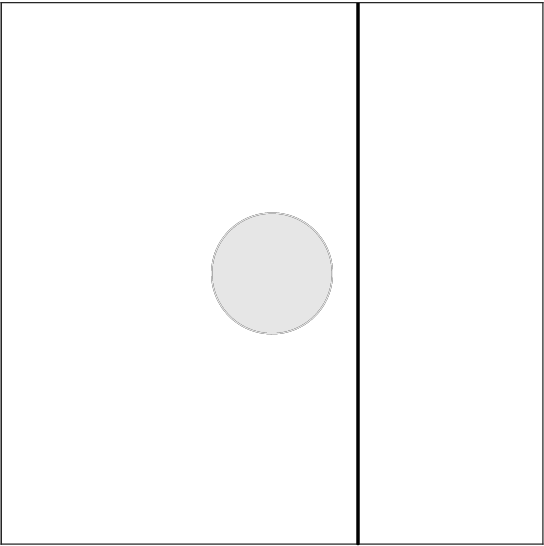}}\hspace{0.3in}
	\subfigure{\includegraphics[width=1.0in]{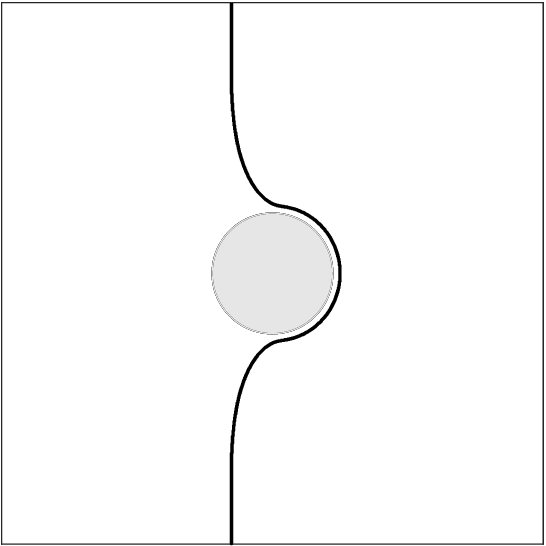}}\hspace{0.3in}
	\subfigure{\includegraphics[width=1.0in]{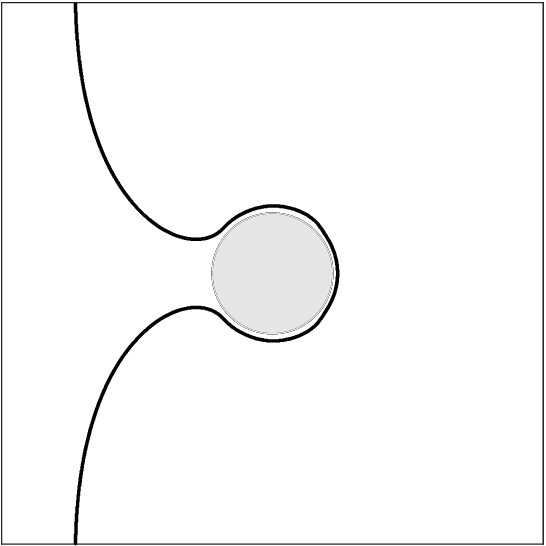}}\hspace{0.3in}
	\subfigure{\includegraphics[width=1.0in]{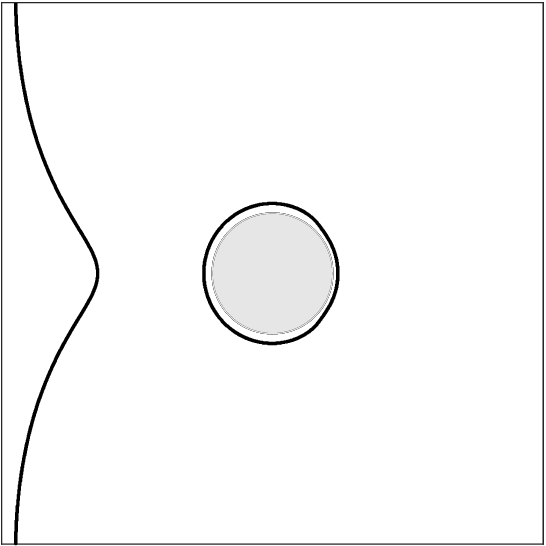}}
	\caption{An edge dislocation bypassing a spherical impenetrable particle under an applied stress by the Orowan loop mechanism. The Burgers vector is in the horizontal direction. The Poisson ratio $\nu=1/3$.  }\label{fig:bypass}
\end{figure}

Simulation result using our threshold dynamics method for dislocation bypassing particle by Orowan loop mechanism is shown in Fig.~\ref{fig:bypass}. In this process, an edge dislocation is driven under an applied stress 
towards an impenetrable spherical particle whose model is given below. As the dislocation approaches the particle, the portion of the dislocation behind the particle is blocked, and the other portions bow forward under the applied stress.  The two dislocation arms on the sides of the particle  continue bow out, and they eventually meet and annihilate. After that, the dislocation pitches off the particle and  leaves behind a dislocation loop around the particle. This is the  Orowan loop bypassing mechanism \cite{anderson2017theory,xiang2004level}.

The spherical particle is modeled by a strong repulsive force acting on any dislocation within the particle and zero on any dislocations outside the particle. We adopt the repulsive force presented in \cite{xiang2004level}:
\begin{equation}
\left\{
	\begin{split}
	&f_0 \hspace{2.15in} \text{if}\hspace{0.1in} r<=R,\\
	&f_0(R+dx -r)^2/dx^2 \hspace{1in} \text{if} \hspace{0.1in}R<r<=R+dx,\\
	&0 \hspace{2.18in} \text{if} \hspace{0.1in}r>R+dx,
	\end{split}	\right .
\end{equation}
where $R$ is the radius of the spherical particle, $r$ is the distance from a point on the dislocation line to the center of the particle, $dx$ is the size of a smooth connecting region. The constant $f_0$ is chosen to be large enough so that the dislocation cannot penetrate the particle. {The radius of the particle $R=0.7\approx 17.5 b$ in the simulation. The numerical grid of the spatial domain is $N=1024$.}

\subsection{Frank-Read source }
We perform simulation for the operation of a Frank-Read source, which is the major mechanism for dislocation multiplication \cite{frank1950multiplication,anderson2017theory}. In this process,  a dislocation segment pinned at both ends bows out under an applied shear stress to generate a series of dislocation loops.

The computational domain is   $[-3\pi,3\pi]\times[-3\pi,3\pi]$, 
discretized into  $N=2048$  grid points in each dimension, i.e. $\Delta x=0.0092$. The dislocation segment with length  $l=2$  is parallel to the $y$ axis, located at  $x_0=\frac{\pi}{3}$ with two end points  $(x_0,y_0)$ and $(x_0,-y_0)$, where $y_0=\frac{l}{2}=1$.
In order to simulate the Frank-Read source, the initial condition $u^0$ is set as
\begin{equation}
	u_0=\left\{
	\begin{aligned}
	 1&  \hspace{0.5in} x_0-\Delta x\leq x\leq x_0, \  -y_0\leq y\leq y_0, \\
	 0&  \hspace{0.5in}\mathrm{otherwise.}	
	\end{aligned}
	 \right.
\end{equation}
This $u^0$ generates a small narrow rectangular counterclockwise dislocation loop. The right vertical segment of the rectangular loop is eliminated in its effect
by adding an extra stress $\sigma_{\rm pin}$ that is generated by a coincident
 dislocation segment with the opposite direction, which serves
to pin the original dislocation segment and cancels its stress field. This pinning stress is \cite{anderson2017theory}:
\begin{equation}
	\sigma_{\rm pin}=-\frac{\mu b}{4\pi(1-\nu)}\left[\frac{(y_0-y)(x-x_0)}{(x-x_0)^2\sqrt{(x-x_0)^2+(y-y_0)^2}}+\frac{(y_0+y)(x-x_0)}{(x-x_0)^2\sqrt{(x-x_0)^2+(y+y_0)^2}}  \right].
\end{equation}
With this pinning stress, the left vertical dislocation segment of the rectangular loop generated by $u_0$ will be operating as a Frank-Read source under an applied stress. We choose the applied stress  $\sigma^{\rm app}=-0.04\mu$ ($-1$ in the dimensionless form).

\begin{figure}[!htb]
	\centering
	\subfigure[]{\includegraphics[width=0.8in]{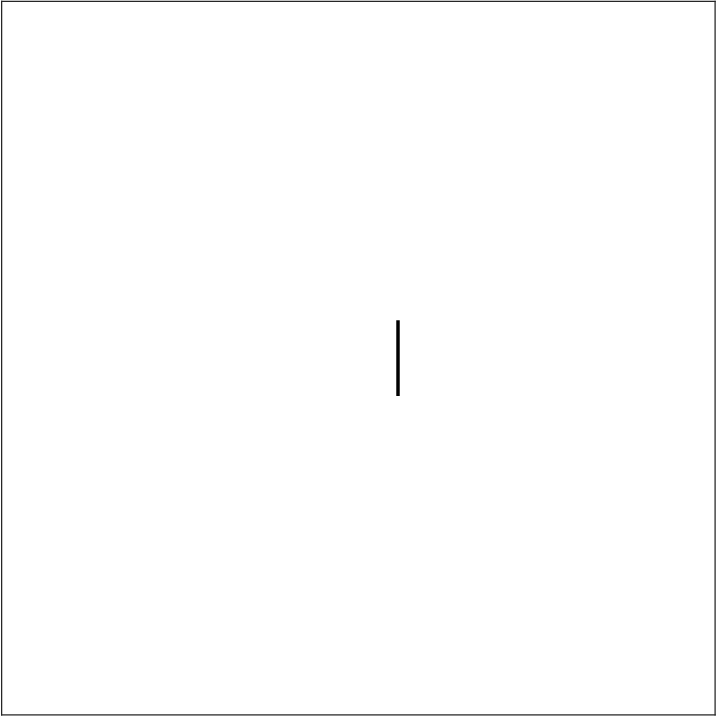}}\hspace{0.25in}
	\subfigure[]{\includegraphics[width=0.8in]{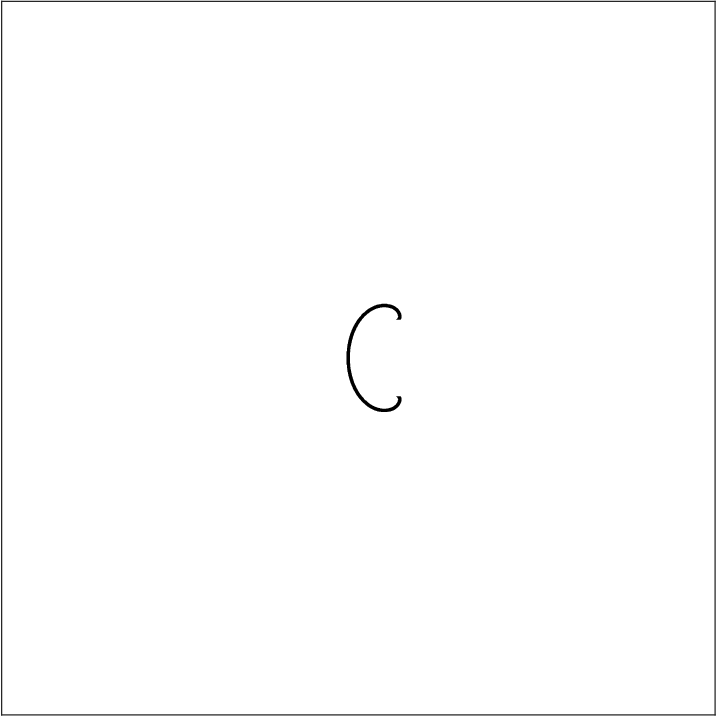}}\hspace{0.25in}
	\subfigure[]{\includegraphics[width=0.8in]{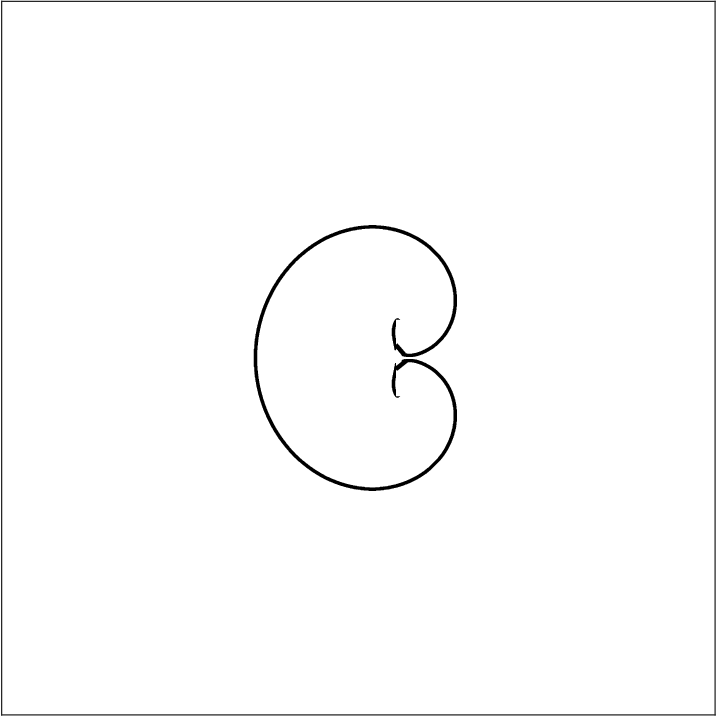}}\hspace{0.25in}
	\subfigure[]{\includegraphics[width=0.8in]{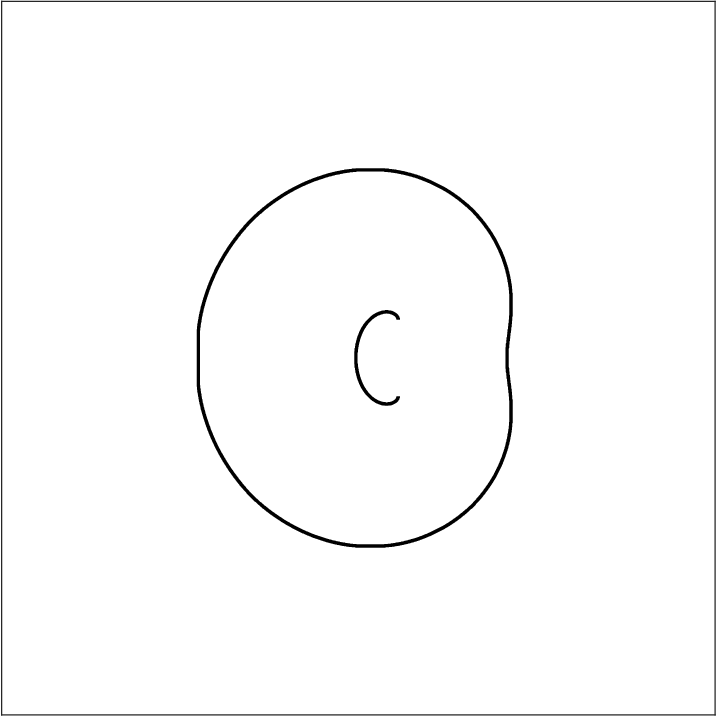}}\hspace{0.25in}
	\subfigure[]{\includegraphics[width=0.8in]{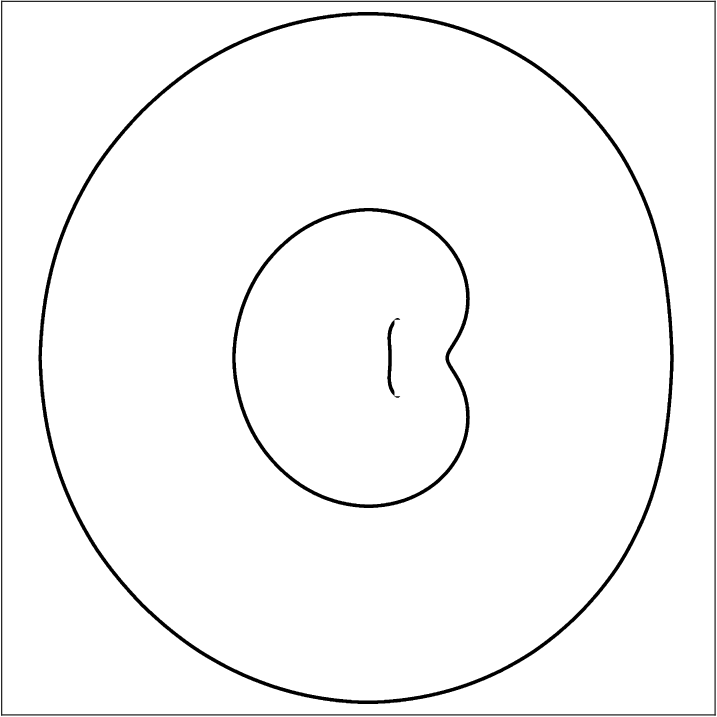}}
	\caption{ Simulation of operation of a Frank-Read source  under an applied stress $\sigma=0.04\mu$. The computational area is  $[-3\pi,3\pi]\times[-3\pi,3\pi]$, and the length of the Frank-Read source is $2$. The Burgers vector is in the horizontal direction. The Poisson ratio $\nu=1/3$. }\label{fig:FRsource}
\end{figure}

Simulation results using our threshold dislocation dynamics method is shown in Fig.~\ref{fig:FRsource}. The dislocation segment bows out to the left under the applied stress, and the two ends of the segment are pinned; see Fig.~\ref{fig:FRsource}(b). When the dislocation segment is heavily bows out, the upper and lower arms meet on the other side of the original segment and the meeting portions are annihilated with each other, and a dislocation loop is pinched off; see Fig.~\ref{fig:FRsource}(c) and (d). After a dislocation loop is pinched off, there is still a dislocation segment inside that connects the two pinned points; see Fig.~\ref{fig:FRsource}(d). This process repeats under the applied stress and more dislocation loops are generated; see  Fig.~\ref{fig:FRsource}(e). This process is the operation of a Frank-Read source~\cite{frank1950multiplication,anderson2017theory}.

\section{Summary}

In this paper, we have developed an efficient threshold dynamics method for dislocation dynamics in a slip plane, in which the spatial operator is essentially an anisotropic fractional Laplacian.  We show that
when setting the time step $\Delta t=\varepsilon$, where $\varepsilon$ is the dislocation core size,
this proposed threshold dislocation dynamics method is able to give correct two leading orders in dislocation velocity,  including both the $O(\log \varepsilon)$ local curvature force and the $O(1)$ nonlocal force due to the long-range stress field generated by the dislocations as well as the force due to the applied stress.
{ This generalizes the result of threshold dynamics formulation with
the kernel of the square root of the Laplacian available in the literature \cite{caffarelli2010convergence}, which is
on the leading order $O(log \Delta t)$ local curvature velocity under the isotropic kernel.}

We have also proposed a numerical method based on stretching of the spatial variables to correct the dislocation mobility
and rescale the dislocation velocity to a larger value for efficient and accurate simulations. This correction method applies generally to any threshold dynamics method for the moving fronts.

We perform numerical simulations using the threshold dislocation dynamics method
for the motion of a straight dislocation under applied stress, shrinking and expanding of dislocation loops, dislocations bypassing particles, and operation of a Frank-Read source.
Simulation results agree with those of theoretic predictions \cite{anderson2017theory} and discrete dislocation dynamics simulations \cite{xiang2003level,xiang2004level}. These simulation results demonstrate that our threshold dislocation dynamics method can indeed correctly capture both the leading order ($O(\log\varepsilon)$) curvature motion and the next order ($O(1)$) long-range interaction for the dynamics of dislocations.

\section*{Acknowledgments}
This work was supported by the Hong Kong Research Grants Council Collaborative Research Fund C1005-19G and the Project of Hetao Shenzhen-HKUST Innovation Cooperation Zone HZQB-KCZYB-2020083. The work of AHWN was also supported by Shenzhen Fund 2021 Basic Research General Programme (project code: JCYJ20210324115400002).

\begin{appendices}

\section{Peierls-Nabarro model for more general cases}\label{sec:PN-general}

Here we give remarks on Peierls-Nabarro model for more general cases, from which more general threshold dislocation dynamics method can be obtained.

{\bf Remark 1.} In a general case, when the Burgers vector is $\mathbf b=(b_1,b_2,0)$, the gradient flow in Eq.~\eqref{eqn:gradient-f} becomes
\begin{equation}\label{eqn:gradient-f1}
\phi_t=-M_p\frac{\delta E}{\delta \phi}=-M_p\left(\sigma_{13}\frac{b_1}{b}+\sigma_{23}\frac{b_2}{b}+\frac{\partial\gamma}{\partial\phi}\right),
\end{equation}
where \cite{xiang2008generalized,zhu-adaptive}
\begin{flalign}
\sigma_{13}^{dis}(x,y)=&\int_\Gamma\left[\frac{\mu b_1}{4\pi(1-\nu)}\frac{(x-\bar{x})}{[(x-\bar{x})^2+(y-\bar{y})^2]^{\frac{3}{2}}} u_x(\bar{x},\bar{y})\right.\nonumber\\
& \ \ \ +\frac{\mu b_1}{4\pi}\frac{(y-\bar{y})}{[(x-\bar{x})^2+(y-\bar{y})^2]^{\frac{3}{2}}} u_y(\bar{x},\bar{y})\nonumber\\
& \ \ \ \left.+\frac{\mu \nu b_2}{4\pi(1-\nu)}\frac{(x-\bar{x})}{[(x-\bar{x})^2+(y-\bar{y})^2]^{\frac{3}{2}}} u_y(\bar{x},\bar{y})\right]d\bar{x}d\bar{y}, \label{eqn:sigma13_2}\\
\sigma_{23}^{dis}(x,y)=&\int_\Gamma\left[\frac{\mu \nu b_1}{4\pi(1-\nu)}\frac{(y-\bar{y})}{[(x-\bar{x})^2+(y-\bar{y})^2]^{\frac{3}{2}}} u_x(\bar{x},\bar{y})\right.\nonumber\\
& \ \ \ +\frac{\mu b_2}{4\pi}\frac{(x-\bar{x})}{[(x-\bar{x})^2+(y-\bar{y})^2]^{\frac{3}{2}}} u_x(\bar{x},\bar{y})\nonumber\\
& \ \ \ \left.+\frac{\mu  b_2}{4\pi(1-\nu)}\frac{(y-\bar{y})}{[(x-\bar{x})^2+(y-\bar{y})^2]^{\frac{3}{2}}} u_y(\bar{x},\bar{y})\right]d\bar{x}d\bar{y},\label{eqn:sigma23_2}
\end{flalign}
and ${\sigma^{\rm app}}=\sigma_{13}^{\rm app}b_1/b+\sigma_{23}^{\rm app}b_2/b$.
Recall that $u$ is the dimensionless disregistry $u=\phi/b$. In this case, the dimensionless evolution of dislocations is still given by Eq.~\eqref{eqn:moelv}, in which $L(u)$ is the dimensionless form of $\sigma_{13}b_1/b+\sigma_{23}b_2/b$.
The Fourier transform of $L(u)$ is
\begin{equation}\label{eqn:sigma13-hat-2}
\widehat{L(u)}=-\frac{1}{2} \left[\left(\frac{1}{1-\nu}\frac{b_1}{b}+\frac{b_2}{b}\right)\frac{k_1^2}{\|\mathbf{k}\|}
+\frac{2b_1b_2\nu}{b^2(1-\nu)}\frac{k_1k_2}{\|\mathbf{k}\|}
+\left(\frac{b_1}{b}+\frac{1}{1-\nu}\frac{b_2}{b}
\right)\frac{k_2^2}{\|\mathbf{k}\|}\right]\hat{u}.
\end{equation}

{\bf Remark 2.} When there are dislocations with multiple Burgers vectors $\mathbf b^{(j)}=(b_1^{(j)},b_2^{(j)},0)$, $j=1,2,\cdots,J$, $J$ disregistry functions $u^{(j)}$s are used, and each $u^{(j)}$ describes dislocations with the Burgers vector $\mathbf b^{(j)}$. The dimensionless evolution is given by
\begin{eqnarray}
u_t^{(j)}=L^{(j)}(u^{(1)},u^{(2)},\cdots,u^{(J)}) -\frac{1}{2\pi \varepsilon}\mathrm{sin}(2\pi u^{(j)})+{\sigma^{app,j}},
\end{eqnarray}
where
\begin{flalign}
L^{(j)}=&\sigma_{13}^{dis}b_1^{(j)}/b^{(j)}+\sigma_{23}^{dis}b_2^{(j)}/b^{(j)},\\
\sigma_{13}^{dis}=&\sum_{j=1}^J \sigma_{13}^{dis,j}, \ \ \
\sigma_{23}^{dis}=\sum_{j=1}^J \sigma_{23}^{dis,j},\\
{\sigma^{app,j}}=&\sigma_{13}^{\rm app}b_1^{(j)}/b^{(j)}+\sigma_{23}^{\rm app}b_2^{(j)}/b^{(j)}.
\end{flalign}
Here $\sigma_{13}^{dis,j}$ and $\sigma_{23}^{dis,j}$ are the stress components generated by dislocations with Burgers vector $\mathbf b^{(j)}$ given by dimensionless form of Eqs.~\eqref{eqn:sigma13_2} and \eqref{eqn:sigma23_2}.

\end{appendices}
	
\bibliographystyle{abbrv}
\bibliography{ref}
\end{document}